# A Directional-ODE Framework for Discretization of Advection-Diffusion Equations

**Amin Jafarimoghaddam**[⋆,1], **Manuel Soler**[†] and **Irene Ortiz**[†]

[⋆,†]*Department of Aerospace Engineering, Universidad Carlos III de Madrid. Avenida de la Universidad, 30, Leganes, 28911 Madrid, Spain* [2]

**Abstract**

We present a novel approach that redefines the traditional interpretation of explicit and implicit discretization methods for solving a general class of advection–diffusion equations (ADEs) featuring nonlinear advection, diffusion operators, and potential source terms. By reformulating the discrete ADEs as directional ordinary differential equations (ODEs) along temporal or spatial dimensions, we derive analytical solutions that lead to novel update formulas. In essence, the information of discrete ADEs is compressed into these directional ODEs, which we refer to as representative ODEs. The analytical update formulas derived from the representative ODEs significantly enhance stability, computational efficiency, and spatiotemporal resolution. Furthermore, we extend the framework to systems with uncertain parameters and coefficients, showcasing its versatility in addressing complex ADEs encountered in modeling and simulation across diverse scientific and engineering disciplines.

## 1  Introduction

Advection–diffusion equations (ADE) are fundamental to understanding diverse processes in science and engineering, such as heat transfer, pollutant dispersion, fluid dynamics, groundwater contaminant transport, chemical reactions in catalytic processes, and nutrient distribution in biological systems. Over the years, researchers have developed various strategies to solve these equations—ranging from analytical and semi-analytical techniques to numerical simulations—to capture the complex dynamics they govern.

---

[1][⋆]Corresponding Author.
[2]Email addresses: ajafarim@pa.uc3m.es (A. Jafarimoghaddam), masolera@ing.uc3m.es (M. Soler), irortiza@ing.uc3m.es (I. Ortiz)

Analytical solutions are generally limited to highly simplified scenarios characterized by assumptions such as linearity [12,17], constant coefficients [3,10], or idealized boundary conditions (see, e.g., [1,15,23,24]). Furthermore, it is important to note that existing analytic series solutions often suffer from ill-conditioning and intractability under certain (typically stiff) conditions (see, e.g., [1,17]). Additionally, some of these solutions fail to reflect realistic or practical scenarios [18].

Semi-analytical methods, such as perturbation, Kernel-based and meshless methods have also been investigated as versatile alternatives for solving a broader class of ADEs. Kernel-based approaches ( [2,7]) are particularly advantageous in handling complex geometries and irregular domains without requiring structured grids. However, they often incur higher computational costs, especially for large-scale or high-dimensional problems. The accuracy and computational efficiency of these methods are typically influenced by the choice of the kernel function, its associated parameters, and the number of support points. In [20], a perturbation-based solution strategy was proposed for general linear ADEs with small diffusion rates. In principle, this approach extends to a broader class of linear ADEs under similar conditions. Nevertheless, the applicability of perturbation-based methods are limited by their reliance on the existence of small parameters, and the number of terms required to achieve a desired solution accuracy is often unknown. Consequently, comprehensive sensitivity analyses—with diverse test cases and parameter settings—are necessary to assess the generalizability of such methods for ADEs.

While analytic, semi-analytic, and meshless (or kernel-based) approaches have advanced our ability to solve advection–diffusion equations under idealized conditions, these methods often struggle when faced with the nonlinearities and uncertainties characteristic of real-world systems. Consequently, in practical applications—where complex boundary conditions, nonlinear behavior, and stochastic parameters prevail—robust numerical methods become indispensable.

Numerical methods, such as finite difference, finite volume, and finite element methods, as well as spectral approaches, are typical tools for solving ADEs in general settings. Given the extensive literature on numerical methods for ADEs, we refer the reader to [6] for a comprehensive review of the various approaches. These methods, however, must navigate challenges such as maintaining computational efficiency, stability, and ensuring solution accuracy, specially under conditions characterized by high Peclet numbers or small diffusion coefficients. Moreover, numerical methods inherently provide a discrete approximation of an evolving field, and thus do not capture the continuous temporal or spatial history between grid points.

In this study, we introduce the directional-ODE discretization approach as a novel paradigm for solving a broad class of advection–diffusion equations (ADEs) that encompass nonlinear advection, nonlinear diffusion operators, and source terms. This method overcomes notable limitations inherent in existing analytical, semi-analytical, and numerical techniques. Unlike traditional analytical and semi-analytical methods—which often require restrictive boundary or initial conditions—the directional-ODE discretization approach accommodates arbitrary conditions while preserving the analytical structure of the solution. Furthermore, whereas conventional numerical methods typically neglect the continuous temporal or spatial history between grid points, our approach integrates this continuous information into the discrete formulation, yielding enhanced accuracy, stability, and computational efficiency.

The directional-ODE discretization approach reinterprets classical discretization schemes by framing them as ODEs in either the temporal or spatial domain. These ODEs, referred to as representative ODEs, serve as reduced-order models that encapsulate the essential structural properties and information of the underlying PDEs. In essence, this paradigm offers a fresh perspective that unifies and redefines various classical discretization frameworks. For clarity and focus, this work emphasizes the finite difference framework, demonstrating how the directional-ODE discretization approach fundamentally reshapes and extends its conventional narrative. Despite its promising implications, this approach has received limited attention in the existing literature, with only a handful of (mostly recent) studies exploring similar ideas, particularly in the context of linear problems.

The foundational contribution in this direction appears in [4], which introduces a semi-implicit finite difference scheme for the one-dimensional linear heat equation. This scheme uniquely treats the midpoint of the second-order derivative term implicitly. Subsequent investigations have demonstrated that the discretization strategy proposed in [4] can be reformulated as an ODE for a class of heat-type equations—including the linear heat equation [13], space- and time-dependent linear diffusion equations [19], linear diffusion-reaction systems [11], and linear Nagumo-type diffusion–reaction equations [14]. In these cases, the resulting representative temporal ODEs are linear and mostly autonomous, admitting simple analytical solutions that yield unconditionally stable forward-time update formulas. Nonetheless, extending this methodology to general ADEs remains a significant open challenge.

Notably, the aforementioned studies represent a special case of the broader directional-ODE discretization approach introduced in this work. More specifically, our framework is not only applicable to a wide class of nonlinear and stochastic ADEs, but it also provides closed-form analytical update formulas for general non-autonomous representative ODEs, with arbitrary order of accuracy, with or without the incorporation of spectral grids. Moreover, unlike the works in [4,11,13,14,19], which consider only the temporal direction, the directional-ODE discretization approach reformulates discrete schemes as ODEs in either the spatial or temporal direction, accompanied by a detailed stability analysis and extensions to problems in arbitrary dimensions. Moreover, for more complex nonlinear representative ODEs, where explicit closed-form solutions are unavailable, we employ the Segmented Adomian Decomposition Method (SADM), a semi-analytical technique that offers a robust solution framework with controllable accuracy.

In summary, the directional-ODE discretization approach enables the derivation of analytic update formulas in both the temporal and spatial directions. Specifically, by *freezing* time, the discrete ADE is transformed into ODEs, which describe the variation of the quantity along the spatial directions. In contrast, by *freezing* space, the discrete ADE is reformulated as ODEs that characterize the temporal evolution of the quantity. To differentiate between these two methodologies, we classify the directional-ODE discretization approach into two categories: **(1)** the temporal-ODE scheme, which typically yields an explicit discretization, and **(2)** the spatial-ODE scheme, which results in an implicit discretization.

The schematic diagrams in Fig. 1 outline the aforementioned methodologies.

Furthermore, Fig. 2 provides a flowchart that encapsulates the core concept of our work, with a particular emphasis on the temporal-ODE scheme.

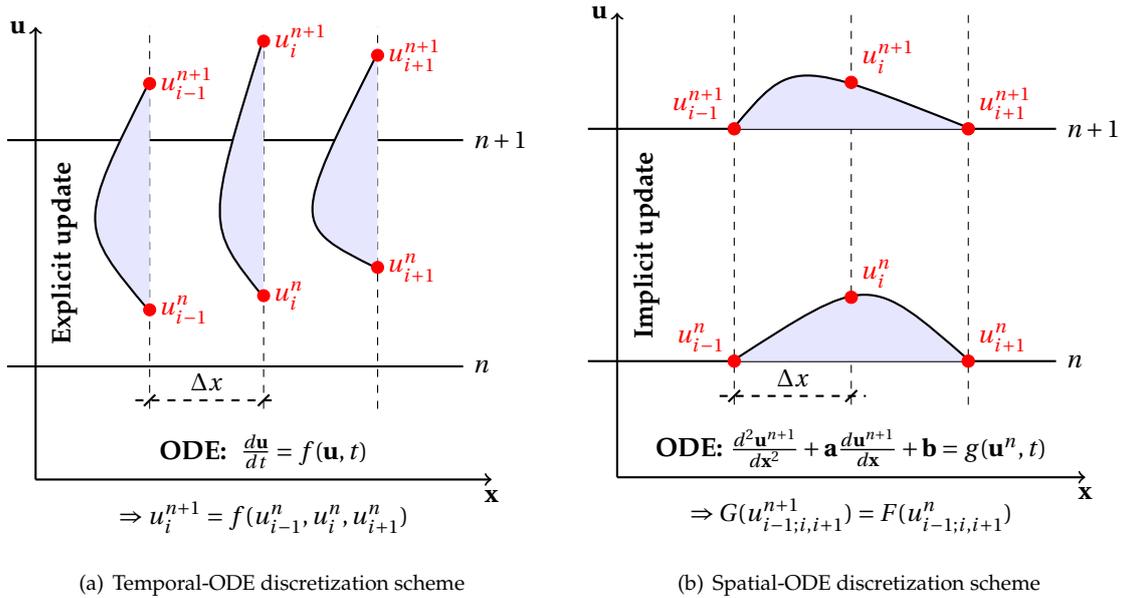

(a) Temporal-ODE discretization scheme  (b) Spatial-ODE discretization scheme

Figure 1: Comparison of temporal-ODE and spatial-ODE schemes

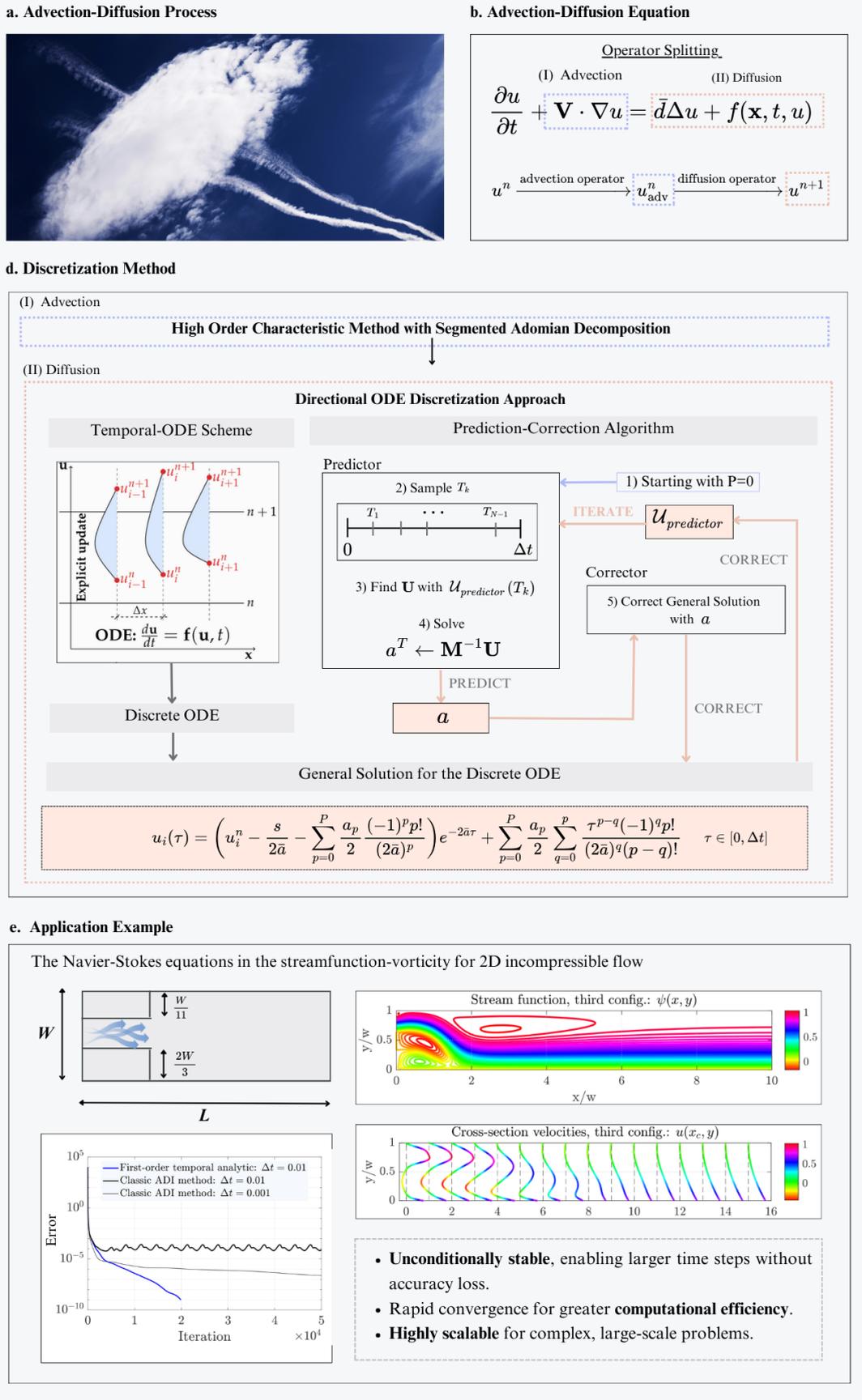

Figure 2: Schematic diagram of the proposed method for solving advection-diffusion equations using the directional-ODE discretization approach presented in this paper

## 2 The Advection-Diffusion Equation

The typical form of a general nonlinear ADE with conserved velocity field, i.e., $\nabla \cdot \mathbf{v} = 0$ is given by:

$$\underbrace{\frac{\partial u(\mathbf{x},t)}{\partial t} + \mathbf{v}(\mathbf{x},t,u) \cdot \nabla u(\mathbf{x},t)}_{\text{advection term}} = \underbrace{\nabla \cdot \left( \bar{d}(\mathbf{x},t,u) \nabla u(\mathbf{x},t) \right)}_{\text{diffusion term}} + \overbrace{f(\mathbf{x},t,u)}^{\text{growth/decay term}}, \qquad (1)$$

where t denotes time and $\mathbf{x}$ denotes the coordinate system in a vector notation, i.e., $\mathbf{x} := (x, y, z)$ in a conventional Cartesian framework, $u$ represents the scalar quantity of interest (such as concentration, temperature or pressure), $\mathbf{v} := (\mathbf{v}_x, \mathbf{v}_y, \mathbf{v}_z)$ represents the velocity field associated with advection of $u$, $\bar{d}$ is the scalar isotropic diffusion coefficient, and $f$ is the system's input/source term which captures external influences or processes that affect the evolution of $\mathbf{u}$.

The advection term in Eq. (1) describes the transport of the quantity driven by $\mathbf{v}$, while the diffusion term models the spread of the quantity due to random particles' motion with the isotropic diffusion $\bar{d}$ controlling rate at which $u$ diffuses through space.

Expanding the diffusion term, the original ADE can be reformulated in the following standard form (see **Supplementary Information, Section 1**):

$$\frac{\partial u}{\partial t} + \mathbf{V} \cdot \nabla u = \bar{d} \Delta u + f. \qquad (2)$$

where, $\mathbf{V} := \mathbf{V}(\mathbf{x}, t, u, \nabla u)$.

To present the directional-ODE discretization approach for general ADEs in a clear and structured manner, we begin by employing an operator splitting technique [5], which not only proves successful across various disciplines of ADEs but also provides a structural perspective through which we can articulate the central concept of the directional-ODE discretization approach in a more systematic manner.

In solving ADEs, the operator splitting method provides a tractable approach by dividing the full system into smaller, more manageable components. Through this technique, the advection and diffusion terms are solved iteratively on each time interval $\Omega_n := [t_n, t_{n+1}], n = 0, 1, ..., N-1$ where $N$ is the total number of temporal grids.

Therefore, on $\Omega_n$, we solve:

$$\begin{cases} \textbf{advection: } \dfrac{\partial u}{\partial t} + \mathbf{V} \cdot \nabla u = 0, \\ \textbf{diffusion: } \dfrac{\partial u}{\partial t} = \bar{d} \Delta u + f. \end{cases} \qquad (3)$$

The solution process from $t_n$ to $t_{n+1}$ can be described as:

$$u^n \xrightarrow{\text{advection operator}} u^n_{\text{adv}} \xrightarrow{\text{diffusion operator}} u^{n+1}. \tag{4}$$

The above scheme aligns with the first-order accurate splitting technique, i.e., $u^{n+1} = \mathcal{D}_{\Delta t}\mathcal{A}_{\Delta t}u^n + O(\Delta t^2)$, where $\mathcal{A}$ and $\mathcal{D}$ are the advection and diffusion operators, respectively.

Higher-order splitting techniques are also available [16, 21, 22]. For example, Strang splitting suggests a second-order accuracy by applying the advection and diffusion operators in an alternating manner, i.e., $u^{n+1} = \mathcal{A}_{\Delta t/2}\mathcal{D}_{\Delta t}\mathcal{A}_{\Delta t/2}u^n + O(\Delta t^3)$.

## 2.1 Revisiting the Nature of the Diffusion Operator

Let us begin by explaining the aforementioned scenarios through an analysis of the diffusion operator involved in the splitting method.

The diffusion equation can be solved using various numerical discretization techniques [6], each with its own advantages and drawbacks. Theoretically, the diffusion operator spreads the initial condition $u^n$ in space exponentially. Nevertheless, existing numerical discretization schemes fail to reflect this structural property in their update formulas. The directional-ODE discretization approach, on the other hand, offers theoretical insights into the diffusion operator by preserving its structural properties. To illustrate this, let us consider the following 1D diffusion equation:

$$\frac{\partial u}{\partial t} = D\frac{\partial^2 u}{\partial x^2} + f(x,t,u). \tag{5}$$

where the constant $D$ is the diffusion coefficient and the initial condition is given as $u(x,0) = u_0(x)$.

### 2.1.1 Temporal-ODE Discretization Scheme

To implement a temporal-ODE scheme, the spatial derivatives are expanded. Unlike classical discretization methods, this approach considers the spatial discrete values as functions of time.

For the 1D diffusion operator, a temporal-ODE scheme can be formed as:

$$\frac{du_i}{dt} = D\frac{u_{i-1}(t) - 2u_i(t) + u_{i+1}(t)}{(\Delta x)^2} + f(x_i, t_n, u_i^n), \quad t \in [t_n, t_{n+1}]. \tag{6}$$

In principle, the neighboring nodes ($u_{i-1}$ and $u_{i+1}$) and the central node ($u_i$) are time-dependent. However, for clarity, we initially fix the neighboring nodes at $t_n$. Therefore, the above ODE can be written as:

$$\frac{du_i}{d\tau} = -Au_i(\tau) + B, \quad u_i(0) = u_i^n, \quad \tau := t - t_n, \quad \tau \in [0, \Delta t], \quad \Delta t := t_{n+1} - t_n. \tag{7}$$

where $A := 2\frac{D}{(\Delta x)^2}$, and $B := D\frac{u_{i-1}^n + u_{i+1}^n}{(\Delta x)^2} + f(x_i, t_n, u_i^n)$.

The analytic solution for the above ODE reads $u_i(\tau) = \frac{B}{A} + (u_i^n - \frac{B}{A})e^{-A\tau}$, giving the explicit analytic update formula at $t_{n+1}$ as:

$$u_i^{n+1} = \frac{B}{A} + \left(u_i^n - \frac{B}{A}\right)e^{-A\Delta t}. \tag{8}$$

### 2.1.2 Spatial-ODE Discretization Scheme

To implement a spatial-ODE scheme, the temporal derivatives are expanded. Unlike classical discretization methods, this approach considers the spatial derivatives as functions of space.

For the 1D diffusion operator, a spatial-ODE scheme can be formed as:

$$\frac{u^{n+1}(x) - u^n(x)}{\Delta t} = D\frac{d^2 u}{dx^2}\bigg|^{n+1} + f(x, t_n, u^n(x)). \tag{9}$$

For clarity, let us assume that $u^n(x)$ is uniformly distributed in space. Consequently, we can express $u^n(x)$ as $u^n(x) := u_i^n$. Similarly, we assume that $u^n(x)$ appearing in $f(x, t_n, u^n(x))$ is uniformly distributed in space. Thus, we can write $f(x, t_n, u^n(x)) := f(x_i, t_n, u_i^n)$.

Therefore, the representative spatial ODE at $t_{n+1}$ can be written as:

$$\frac{d^2 u}{dx^2} - Au = B. \tag{10}$$

where, $A = \frac{1}{D\Delta t}$ and $B = -\frac{u_i^n}{D\Delta t} - \frac{f(x_i, t_n)}{D}$.

In order to form a spatial-ODE scheme, we solve the above ODE on $\Omega := [x_i - \Delta x, x_i + \Delta x]$. In other words, the neighboring nodes are selected as the required boundary conditions:

$$u(x_i - \Delta x) = u_{i-1}^{n+1}, \quad u(x_i + \Delta x) = u_{i+1}^{n+1}. \tag{11}$$

The analytic solution for the above ODE reads $u(x) = -\frac{B}{A} + c_1 e^{\sqrt{A}x} + c_2 e^{-\sqrt{A}x}$, together with the boundary conditions, giving an implicit update formula at $u(x_i) \equiv u_i^{n+1}$ as:

$$u_i^{n+1} = -\frac{B_i}{A} + \frac{u_{i-1}^{n+1} + u_{i+1}^{n+1} + 2\frac{B_i}{A}}{e^{\sqrt{A}\Delta x} + e^{-\sqrt{A}\Delta x}}. \tag{12}$$

Notably, we can rearrange the above analytic discrete scheme into a triangular form.

**Remark:** For higher-dimensional discretization, the spatial-ODE scheme can be implemented using the alternating directional implicit (ADI) method.

Finally, it is evident that the exponential property of the diffusion operator is embedded in both temporal-ODE and spatial-ODE schemes.

## 3 The General Solution Approach

With a focus on generality, this paper employs the analytic characteristic method to handle the advection operator, while the directional-ODE discretization approach is exclusively applied to the diffusion operator. Notably, the directional-ODE discretization approach can be adapted to handle the advection operator and also the entire ADE without any splitting approach (see **Supplementary Information, Section 6**).

In the following, we begin by applying the method of characteristics to a general nonlinear advection operator. Next, we propose a general multi-stage, high-order, closed-form temporal-ODE discretization scheme for the diffusion operator. While this multi-stage approach is applicable to diffusion operators with nonlinear diffusion coefficients and nonlinear source terms, we also introduce a semi-analytic method for handling nonlinear diffusion operators, which we show can reduce computational time in nonlinear scenarios with a competitive level of accuracy.

For clarity, the proposed methodology is explained in 1D, and the 3D version is referred to in **Supplementary Information, Section 2**.

### 3.1 Advection Operator: The Method of Characteristics

To solve the advection equation on the interval $[t_n, t_{n+\frac{1}{2}}]$, we employ the method of characteristics.

The advection equation in 1D can be written as:

$$\frac{\partial u}{\partial t} + \mathbf{V}_x \frac{\partial u}{\partial x} = 0, \quad u(x, t_n) = u^n. \tag{13}$$

The method of characteristics requires solving the following ODE to trace the particle paths:

$$\frac{dx}{dt} = \mathbf{V}_x, \quad x(t_n) = x_0. \tag{14}$$

Notably, the solution remains constant along the characteristic path:

$$u(x, t) = u(x_0, t_n). \tag{15}$$

where $u(x_0, t_n)$ is the initial condition.

For a general semi-analytic solution, one can apply methods such as the SADM [8, 9]. On using the SADM, the characteristic solution is obtained as:

$$x_0 = x - \mathcal{L}^{-1}[\mathbf{V}_x] := x - \int_{t_n}^{t_{n+\frac{1}{2}}} \mathcal{A}_0 \, dt - \int_{t_n}^{t_{n+\frac{1}{2}}} \mathcal{A}_1 \, dt - \ldots \quad (16)$$

where $\mathcal{A}_j$ represents Adomian terms (see Section 3.3.1).

Notably, depending on the desired accuracy, one may also adopt Runge-Kutta methods or the Euler approximation as the simplest approach to solving the above ODE on the interval $[t_n, t_{n+\frac{1}{2}}]$, resulting in: $x_0 = x - \mathbf{V}_x \frac{\Delta t}{2}$.

Therefore, the general solution to the advection operator involves the following coordinate mapping at each time step:

$$\mathcal{A}_{\frac{\Delta t}{2}} u^n := u\left(x - \mathcal{L}^{-1}[\mathbf{V}_x]\right). \quad (17)$$

The coordinate mapping described above can be achieved using various interpolation methods (such as linear, spline, and cubic interpolations), which are readily available in programming languages such as MATLAB.

## 3.2 Diffusion Operator: Multi-Stage Temporal-ODE Discretization Scheme

In this section, we present a general multi-stage temporal-ODE scheme for the diffusion operator.

For the 1D diffusion operator, a temporal-ODE scheme can be formed as:

$$\frac{du_i}{d\tau} = D \frac{u_{i-1}(\tau) - 2u_i(\tau) + u_{i+1}(\tau)}{(\Delta x)^2} + f(x_i, t_n, u_i^n), \quad \tau \in [0, \Delta t]. \quad (18)$$

Since the neighboring nodes are functions of time, we assume: $u_{i-1}(\tau) + u_{i+1}(\tau) =: \mathcal{U} = \sum_{p=0}^{P} a_p \tau^p$,

Defining $\bar{a} := \frac{D}{(\Delta x)^2}$, the representative ODE becomes:

$$\frac{du}{d\tau} = \bar{a} \Big( \sum_{p=0}^{P} a_p \tau^p \Big) - 2\bar{a} u + s. \quad (19)$$

The closed-form solution reads:

$$u_i(\tau) = \left( u_i^n - \frac{s}{2\bar{a}} - \sum_{p=0}^{P} \frac{a_p}{2} \frac{(-1)^p p!}{(2\bar{a})^p} \right) e^{-2\bar{a}\tau} + \sum_{p=0}^{P} \frac{a_p}{2} \sum_{q=0}^{p} \frac{\tau^{p-q}(-1)^q p!}{(2\bar{a})^q (p-q)!} + \frac{s}{2\bar{a}}. \quad (20)$$

where, $s := f(x_i, t_n, u_i^n)$ is known at $t_n$.

In above, $a_p$ are calculated according to an intermediate sampling over $[0, \Delta t]$, leading to a *Predictor-Corrector Algorithm*.

**Predictor Step:** We initially predict intermediate values of $\mathcal{U}(\tau)$ using the update formula associated with $P = 0$:

$$u_i^{predictor}(\tau) =: \mathcal{F}_{predictor}(u_i^n, \tau) = \left(u_i^n - \frac{s}{2\bar{a}} - \frac{a_0}{2}\right)e^{-2\bar{a}\tau} + \frac{a_0}{2} + \frac{s}{2\bar{a}}. \tag{21}$$

Therefore, the sample points are:

$$\mathcal{S} := \{(0, \mathcal{U}_{predictor}^n), (T_1, \mathcal{U}_{predictor}^{n+\frac{1}{N}}), ..., (T_{N-1}, \mathcal{U}_{predictor}^{n+\frac{(N-1)}{N}}), (T_N, \mathcal{U}_{predictor}^{n+1})\}. \tag{22}$$

Here, $N \in \mathbb{N}$ denotes the order of the required polynomial approximation for the neighboring nodes. Moreover, $T_k$, for $k = 1, ..., N$, represents the intermediate time steps at which the samples are collected. For example, uniform time stepping is given by $T_k = k\frac{\Delta t}{N}$, $k = 1, ..., N$. For increased accuracy, non-uniform time stepping, such as Chebyshev nodes, can be employed, giving $T_k = \frac{\Delta t}{2}\left(1 - \cos\left(\frac{k\pi}{N}\right)\right)$, $k = 1, ..., N$.

**Corrector Step:** Using $N+1$ sample points, we can symbolically solve the following system to determine the polynomial coefficients $a_p$ for $p = 0, 1, ..., N = P$, in terms of the predicted neighboring nodes:

$$\begin{pmatrix} 1 & 0 & 0 & \cdots & 0 \\ 1 & T_1^1 & T_1^2 & \cdots & T_1^P \\ 1 & T_2^1 & T_2^2 & \cdots & T_2^P \\ \vdots & \vdots & \vdots & \ddots & \vdots \\ 1 & T_{N-1}^1 & T_{N-1}^2 & \cdots & T_{N-1}^P \\ 1 & T_N^1 & T_N^2 & \cdots & T_N^P \end{pmatrix} \begin{pmatrix} a_0 \\ a_1 \\ a_2 \\ \vdots \\ a_{P-1} \\ a_P \end{pmatrix} = \begin{pmatrix} \mathcal{U}_{predictor}^n \\ \mathcal{U}_{predictor}^{n+\frac{1}{N}} \\ \mathcal{U}_{predictor}^{n+\frac{2}{N}} \\ \vdots \\ \mathcal{U}_{predictor}^{n+\frac{N-1}{N}} \\ \mathcal{U}_{predictor}^{n+1} \end{pmatrix}. \tag{23}$$

Next, we exploit $a_p, p = 0, 1, ..., P$ into the corresponding $P^{th}$-order closed-form solution to compute $u_i^{n+1}$.

Notably, the corrector step can be interpreted as additional predictor steps. After completing the initial predictor step and computing the polynomial coefficients $a_p, p = 0, 1, ..., P$, the corresponding $P^{th}$-order closed-form solution can undergo the same sampling process as the initial step. In other words, the corrector step can be embedded in a loop to iteratively refine the polynomial coefficients $a_p, p = 0, 1, ..., P$. Physically, this process resembles allowing the diffusion operator (i.e., the heat equation) sufficient time to bring neighboring nodes to an equilibrium state at each time step.

## 3.3 Temporal-ODE Discretization Scheme with Nonlinear Diffusion Coefficient

In many practical scenarios, the diffusion coefficient is nonlinear, i.e., $D := D(u)$. The general multi-stage temporal-ODE scheme (outlined in the previous section) is also applicable to cases with a nonlinear diffusion coefficient by considering $D := D(u_i^n)$. However, a more accurate approach is to solve the following representative nonlinear ODE, accounting for the variation of the diffusion coefficient over $[t_n, t_{n+1}]$:

$$\frac{du_i}{d\tau} = D(u_i(\tau))\frac{u_{i-1}(\tau) - 2u_i(\tau) + u_{i+1}(\tau)}{(\Delta x)^2} + f(x_i, t_n, u_i^n), \quad \tau \in [0, \Delta t]. \qquad (24)$$

Following a zeroth-order approximation for the neighboring nodes, the above ODE can be written as:

$$\frac{du_i}{d\tau} = D(u_i(\tau))(Au_i(\tau) + B) + C. \qquad (25)$$

where $A = -\frac{2}{(\Delta x)^2}$, $B = \frac{1}{(\Delta x)^2}(u_{i-1}^n + u_{i+1}^n)$, and $C = f(x_i, t^n, u_i^n)$.

Closed-form solutions can be derived for specific forms of $D(u_i(\tau))$. However, a more general approach involves utilizing a semi-analytic method that can handle any type of nonlinear equation. To this end, we employ the SADM, which is detailed in the following section.

### 3.3.1 Segmented Adomian Decomposition Method (SADM)

Briefly, the SADM is a modified version of ADM which resumes the decomposition over each time interval $[t_n, t_{n+1}]$ (see [8,9] for detailed information).

Following the SADM, we initially write the ODE in the following form $\mathscr{L}[u(t)] = \mathscr{N}[u(t)]$, where $\mathscr{L}$ denotes the linear operator, such as $\frac{d}{dt}$, and $\mathscr{N}$ denotes the nonlinear operator, typically representing the right-hand side (RHS) of the equation.

Next, the solution $u(t)$ and the nonlinearity $\mathscr{N}[u(t)]$ are decomposed over $[t_n, t_{n+1}]$ as:

$$u(t) = u_0(t) + u_1(t) + u_2(t) + ... = \sum_{n=0}^{\infty} u_n(t), \quad \mathscr{N}[u(t)] = \sum_{n=0}^{\infty} \mathscr{A}_n(u(t)). \qquad (26)$$

where, $\mathscr{A}_n(u_0(t), u_1(t), ..., u_n(t)) = \frac{1}{n!}\frac{d^n}{dp^n}\left[\mathscr{N}\left(\sum_{n=0}^{\infty} u_n p^n\right)\right]_{p=0}$, with $u_n(t) = \int_{t_n}^{t} \mathscr{A}_{n-1}\,dt$, and $\mathscr{A}_n$ are the *Adomian* terms.

Notably, the SADM has been shown to be highly stable across various nonlinear equations.

## 3.4 Stability Analysis of the Diffusion Operator

In this section, we briefly examine the conditions under which a general nonlinear representative ODE for the diffusion operator remains stable. Stability is evaluated by analyzing the stability of the stationary point, which corresponds to the solution as $t \to \infty$.

Physically, the source term in the ADE plays a fundamentally different role from the diffusion operator. Growth terms, particularly unbounded or unconstrained functions, often cause the system to exhibit unbounded behavior, leading to unstable stationary points when they exist. As a result, the solution typically does not converge to a constant value as $t \to \infty$. In the presence of such growth terms, steady-state solutions are generally absent, and stability analysis becomes less relevant, shifting the focus instead to accuracy analysis. In contrast, decay terms generally act as stabilizers, facilitating the convergence of the solution to a constant value as $t \to \infty$[3].

Therefore, without loss of generality, by focusing solely on the diffusion operator, we provide a more rigorous stability framework without the complications introduced by growth sources.

The stationary point for Eq. (25) – ignoring the source term – can be obtained by solving the following algebraic equation:

$$D(u^*)\left(Au^* + B\right) =: F(u^*) = 0, \quad \Rightarrow u^* := -\frac{B}{A}. \tag{27}$$

The stability of the stationary point can be analyzed using linear perturbation theory. To this end, by considering a small perturbation $u = u^* + \epsilon(t)$, the perturbed ODE and the associated solution become:

$$\frac{du^*}{d\tau} + \frac{d\epsilon}{d\tau} = F(u^* + \epsilon) = F(u^*) + \epsilon \frac{dF(u^*)}{du^*} + o(\epsilon^2), \quad \Rightarrow \frac{d\epsilon}{d\tau} = \epsilon \frac{dF}{du^*}, \quad \Rightarrow \epsilon(t) = c_1 e^{\frac{dF}{du^*}\tau}. \tag{28}$$

The stationary point is stable if $\frac{dF}{du^*} < 0$. Therefore, we write $\frac{dF}{du^*} = \frac{dD(u^*)}{du^*}\left(Au^* + B\right) + AD(u^*)$. Since $u^* := -\frac{B}{A}$, it follows that $\frac{dF}{du^*} = AD(u^*)$. Because the diffusion term is assumed positive and $A < 0$, we have $\frac{dF}{du^*} < 0$, implying that the stationary point is stable.

Notably, since the directional-ODE discretization approach represents update formulas over discrete intervals $[t_n, t_{n+1}]$, the condition $t \to \infty$ is equivalent to $\Delta t \to \infty$, corresponding to the steady-state solution. In other words, in a discrete format, we check whether the discrete scheme can converge to a fixed value for arbitrary time steps. Therefore, the stability of a discrete format refers to its asymptotic behavior as $\Delta t \to \infty$.

---

[3]From another perspective, the split equation $\frac{\partial u}{\partial t} = D\frac{\partial^2 u}{\partial x^2} + f(x,t,u)$ can be decomposed into $\frac{\partial u}{\partial t} = D\frac{\partial^2 u}{\partial x^2}$ and $\frac{\partial u}{\partial t} = f(x,t,u)$. The overall stability of the system relies on the stability behavior of both sub-equations. As discussed, the decay component is typically associated with stable stationary points, while growth terms generally preclude steady-state solutions. Consequently, in such cases, the emphasis should shift from stability to accuracy analysis.

In **Supplementary Information, Section 3**, we present an exclusive stability analysis of the multi-stage temporal-ODE scheme, demonstrating that the discrete approach is unconditionally stable.

## 3.5 Temporal-ODE Discretization Scheme with Uncertain/Stochastic Diffusion Coefficient

In this section, we examine the diffusion operator and show that the temporal-ODE approach allows for the direct integration of the probabilistic dynamics of the uncertain coefficient into the discrete formulation. Because this formulation maintains the structural properties of the diffusion operator, it is expected to produce more accurate and realistic results. Additionally, this approach can be easily extended to higher dimensions and adapted to various directional-ODE schemes. A detailed analysis is provided in **Supplementary Information, Section 4** for the $P^{th}$-order probabilistic solution accounting for the temporal-ODE scheme. Here, however, we present the simplest probabilistic discrete update formula for the most basic case, which is the zeroth-order solution in the absence of a source term. The solution is given by:

$$\mathbb{E}[u_i^{n+1}] = \left(u_i^n - \frac{a_0}{2}\right)\frac{e^{-2\tilde{D}_l \bar{\bar{a}} \Delta t} - e^{-2\tilde{D}_u \bar{\bar{a}} \Delta t}}{2\bar{\bar{a}}\Delta t(\tilde{D}_u - \tilde{D}_l)} + \frac{a_0}{2}. \tag{29}$$

where $\bar{\bar{a}} := \frac{1}{(\Delta x)^2}$, and $a_0 := u_{i-1}^n + u_{i+1}^n$, and $\tilde{D}$ represents the uncertain diffusion coefficient. For simplicity, we have assumed that $\tilde{D}$ has a uniform probability density function (PDF) over the range $\tilde{D} \in [\tilde{D}_l, \tilde{D}_u]$, where the subscripts $l$ and $u$ represent the lower and upper bounds, respectively.

# 4 Case Studies

This section is structured to highlight various features of the proposed directional-ODE schemes.

In **Example 1**, we address the 1D Burgers' equation, comparing the performance of several directional-ODE schemes with the available closed-form solution and the classic implicit discretization scheme. Particularly, **Example 1** demonstrates that both spatial-ODE and temporal-ODE schemes can outperform classical implicit methods.

In **Example 2**, we solve the 2D propagation of particles under smooth and noisy velocity fields with a nonlinear diffusion coefficient and growth term. Here, the performance of the zeroth-order linear directional-ODE approach is compared to a nonlinear directional-ODE approach using the SADM. Notably, **Example 2** highlights a nonlinear scenario where the representative ODE does not provide a practical closed-form solution.

Finally, in **Example 3**, we solve the 2D Navier-Stokes equations for various back-step flow configurations and compare the performance of the directional-ODE approach with the classic alternating directional implicit (ADI) method. Particularly, **Example 3** demonstrates not only the superiority of the directional-ODE schemes over classical implicit methods but also their effectiveness in handling highly nonlinear situations.

## 4.1 Example 1: 1D Burgers' Equation

Our first test case includes 1D Burgers' equation with specific boundary and initial conditions as described below:

$$\frac{\partial u}{\partial t} + u \frac{\partial u}{\partial x} = \nu \frac{\partial^2 u}{\partial x^2}, \quad |x| \le 1, \quad t > 0, \quad u(x,0) = -\sin(\pi x), \quad u(1,t) = u(-1,t) = 0. \tag{30}$$

where $\nu$ is the dynamic viscosity.

The analytic solution for the above equation with specified boundary conditions is [17]:

$$u(x,t) = -\frac{\int_{-\infty}^{\infty} \sin(\pi(x-\eta)) f(x-\eta) e^{\frac{-\eta^2}{4\nu t}} d\eta}{\int_{-\infty}^{\infty} f(x-\eta) e^{\frac{-\eta^2}{4\nu t}} d\eta}, \quad f(y) = e^{\frac{\cos(\pi y)}{2\pi \nu}}. \tag{31}$$

To isolate the performance of discrete schemes for the diffusion operator, a first-order operator splitting method was adopted. Specifically, the advection operator was solved over the interval $\Omega_n = [t_n, t_{n+1}], n = 0, \ldots, N-1$, exploiting the method of characteristics, coupled with the diffusion operator. However, the diffusion operator was handled using the following discretization schemes: the classic implicit method[4], spatial-ODE scheme, zeroth-order temporal-ODE scheme, and first-order temporal-ODE scheme (with or without loop correction). Notably, for the case of loop correction, we used 20 iterations per corrector step. Moreover, for the analytic formula, the domain in the $\eta$-direction was set to 1000 and discretized using $10^6$ nodes, while the $x$-direction was discretized over 1000 nodes. The dynamic viscosity, $\nu$, was set to a constant value of $\nu = 0.005$ for all the cases studied.

The averaged error was computed over 10 time steps using the following formula:

$$\bar{E} = \frac{1}{10 N_x} \sum_{n=1}^{10} \sum_{i=1}^{N_x} \left( u_{\text{analytic},i}^n - u_i^n \right)^2,$$

where $N_x$ represents the number of discrete nodes in the $x$-direction. By varying $N_x$, comparisons were performed over a wide range of $\lambda := D \frac{\Delta t}{(\Delta x)^2}$, which represents the well-known CFL condition.

---

[4]In the classic implicit scheme, the diffusion operator is discretized as: $\frac{u_i^{n+1} - u_i^n}{\Delta t} = \nu \frac{u_{i-1}^{n+1} - 2u_i^{n+1} + u_{i+1}^{n+1}}{(\Delta x)^2}$.

Fig. 3(a) shows the variation of the averaged error for the studied discrete schemes. From this figure, it is evident that the spatial-ODE scheme gives better results compared to the classic schemes for lower values of $\lambda$, which links to its enhanced performance with fewer spatial nodes $N_x$. This behavior is expected since, as discussed, the spatial-ODE scheme directly solves the representative ODE in the spatial direction. Moreover, from this figure, we observe that zeroth-order and first-order temporal-ODE schemes (without loop correction) exhibit better performance for lower values of $\lambda$. However, while stable, for larger values of $\lambda$ they show lower accuracy compared to other studied methods. In contrast, the accuracy associated with the first-order temporal-ODE scheme improves for larger values of $\lambda$, specifically marking the temporal-ODE schemes with loop correction as a remarkable choice for achieving higher accuracy levels.

In addition, for the steady-state Burgers' equation, which admits a closed-form solution, we solve for the uncertain dynamic viscosity $v$ with a uniform probability distribution $v \in [0.1, 0.9]$ using the discrete formula presented in Eq. (29) and compare it with a direct sampling approach. To quantify the discrepancy, we define the average error as the absolute difference between the expected solution obtained from sampling and that of Eq. (29), averaged over the spatial domain: Average Error $:= \frac{1}{L} \int_0^L \left| \mathbb{E}[u_{\text{sample}}(x)] - u_{\text{exact}}(x) \right| dx.$

where $\mathbb{E}[u_{\text{sample}}(x)]$ represents the mean of the sampled solutions at each spatial location $x$, $u_{\text{exact}}(x)$ denotes the analytic solution, and $L$ is the length of the spatial domain.

As observed in Fig. 3(b), the expected solution obtained from the direct sampling method converges to the solution given by the probabilistic update formula as the number of samples increases.

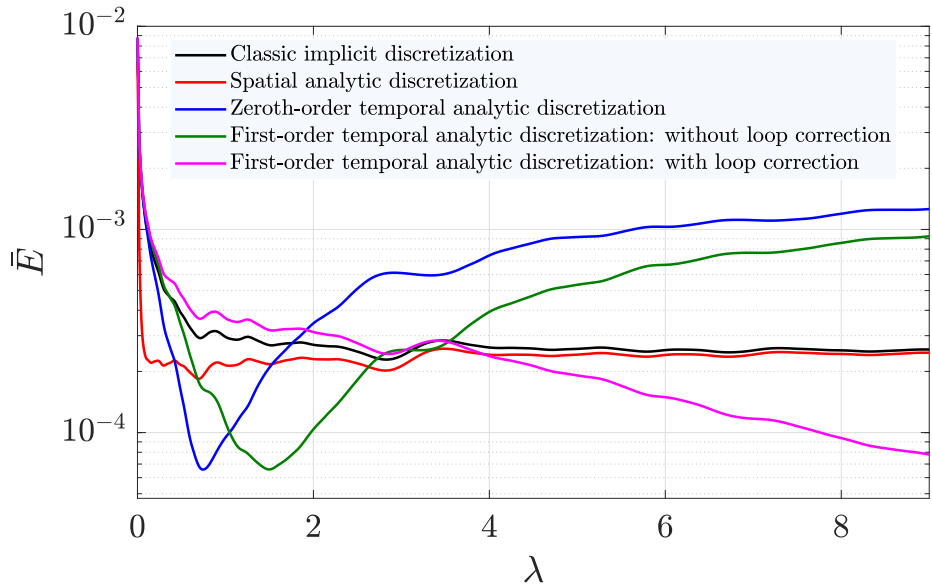

(a) Performance comparison of different discretization schemes for the 1D Burgers' equation.

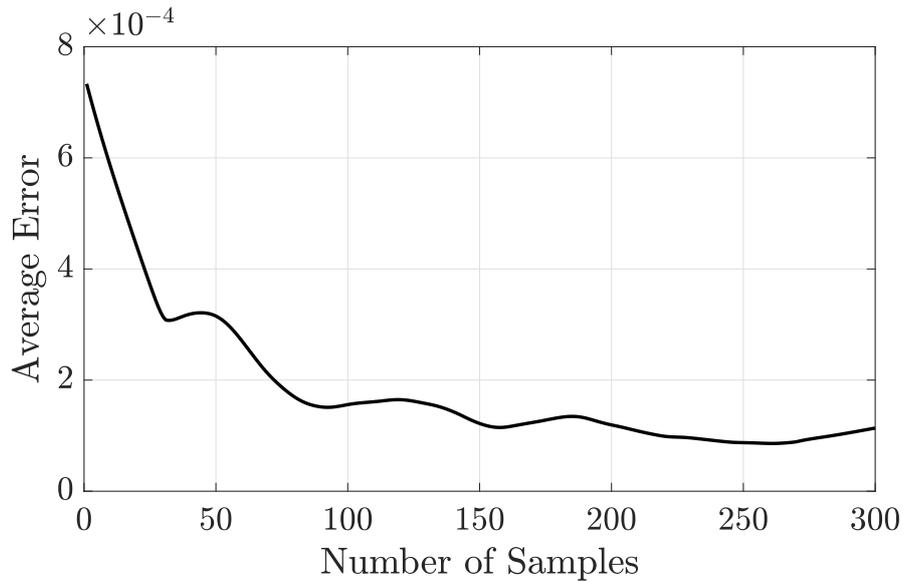

(b) Comparison of the probabilistic solution using the zeroth-order update formula (Eq. 29) and the direct sampling approach.

Figure 3: **Example 1** — Performance and probabilistic solution comparisons for the 1D Burgers' equation.

## 4.2 Example 2: 2D Propagation of Particles

The second test case involves the 2D propagation of a linear concentration field with a nonlinear diffusion and growth term. The considered ADE is:

$$\frac{\partial u}{\partial t} + \mathbf{v} \cdot \nabla u = \nabla(D(u)\nabla u) + su. \tag{32}$$

where $D(u) = \frac{D_0}{1+\beta u}$ with $D_0 = 0.001$, and $\beta > 0$. Moreover, the growth constant is set to $s = 0.01$.

We choose the following wind field ensuring conservation of the flow continuity equation:

$$\mathbf{v}_x = 1.5(1 + kr_{y,1})\sin(y) + (1 + kr_{y,2})\cos(y),$$

$$\mathbf{v}_y = 1.5(1 + kr_{x,1})\sin(x) + (1 + kr_{x,2})\cos(x).$$

where $r_{x,1}$, $r_{x,2}$, $r_{y,1}$, and $r_{y,2}$ represent random parameters drawn from a normal distribution. Additionally, $k$ is set to either $k = 0$ (smooth velocity field) or $k = 1$ (highly noisy velocity field).

Following a double-step splitting method, over $[t_n, t_{n+1}]$, we solve:

$$\textbf{Advection:} \ \frac{\partial u}{\partial t} + \mathbf{V}_x \frac{\partial u}{\partial x} + \mathbf{V}_y \frac{\partial u}{\partial y} = 0, \quad \textbf{Diffusion:} \ \frac{\partial u}{\partial t} = D(u)\left(\frac{\partial^2 u}{\partial x^2} + \frac{\partial^2 u}{\partial y^2}\right) + su. \tag{33}$$

where $\mathbf{V}_x := \mathbf{v}_x - \frac{\partial D}{\partial u}\frac{\partial u}{\partial x}$, and $\mathbf{V}_y := \mathbf{v}_y - \frac{\partial D}{\partial u}\frac{\partial u}{\partial y}$.

The second-order Strang splitting method is adopted. The advection operator is solved using the method of characteristics, while the diffusion operator is addressed through various temporal-ODE schemes. Specifically, we first solve the diffusion operator using a nonlinear temporal-ODE scheme implemented by means of the SADM. The results are then compared against those obtained using the zeroth-order temporal-ODE scheme. For this comparison, a second-order temporal-ODE scheme with loop correction and a highly refined temporal grid is used as the reference solution.

Since it can be shown that the closed-form solution for the representative ODE is an implicit formula, semi-analytic methods such as the SADM can be adopted alternatively. Following the SADM, we obtain:

$$\begin{aligned}
u_0 &= u(0) := u_{i,j}^n, \\
\mathscr{A}_0 &= D_0 \frac{Au_0 + B}{1 + \beta u_0} + C \Rightarrow u_1(\tau) = \int_0^\tau \mathscr{A}_0 d\tau, \\
\mathscr{A}_1 &= \frac{D_0 u_1(\tau)(A - \beta B)}{(\beta u_0 + 1)^2} \Rightarrow u_2(\tau) = \int_0^\tau \mathscr{A}_1 d\tau, \\
&\ldots
\end{aligned} \tag{34}$$

Therefore, $u(\tau) = u_0 + u_1(\tau) + u_2(\tau) + \ldots$, and the update formula (i.e., $u^{n+1}$) is obtained at $\tau = \Delta t$.

The initial particles' concentration is represented by multiple intersecting finite lines with a uniform concentration of 0.01 and a small thickness of 0.02. Moreover, the boundary conditions are set to zero.

Fig. 4, and 5 show the evolution of particles' concentration for the case of $\beta = 10$ with smooth and noisy velocity fields respectively. Notably, the graphed solutions were obtained using the first-order temporal-ODE scheme with loop correction. Particularly, Fig. 5 demonstrates the stability of the solver with a highly noisy velocity field ($k = 1$). Fig. 6 presents a comparison of errors between two temporal-ODE schemes: the zeroth-order closed-form scheme, which approximates the nonlinear diffusion at $t_n$ and reduces the problem to a linear ODE, and the third-order semi-analytic scheme, which accounts for time-dependent nonlinear diffusion over $[t_n, t_{n+1}]$. For this analysis, a highly refined second-order closed-form linear temporal-ODE scheme with loop correction was used as the reference solution. Errors were calculated as the absolute difference across all spatial grid points at various time steps.

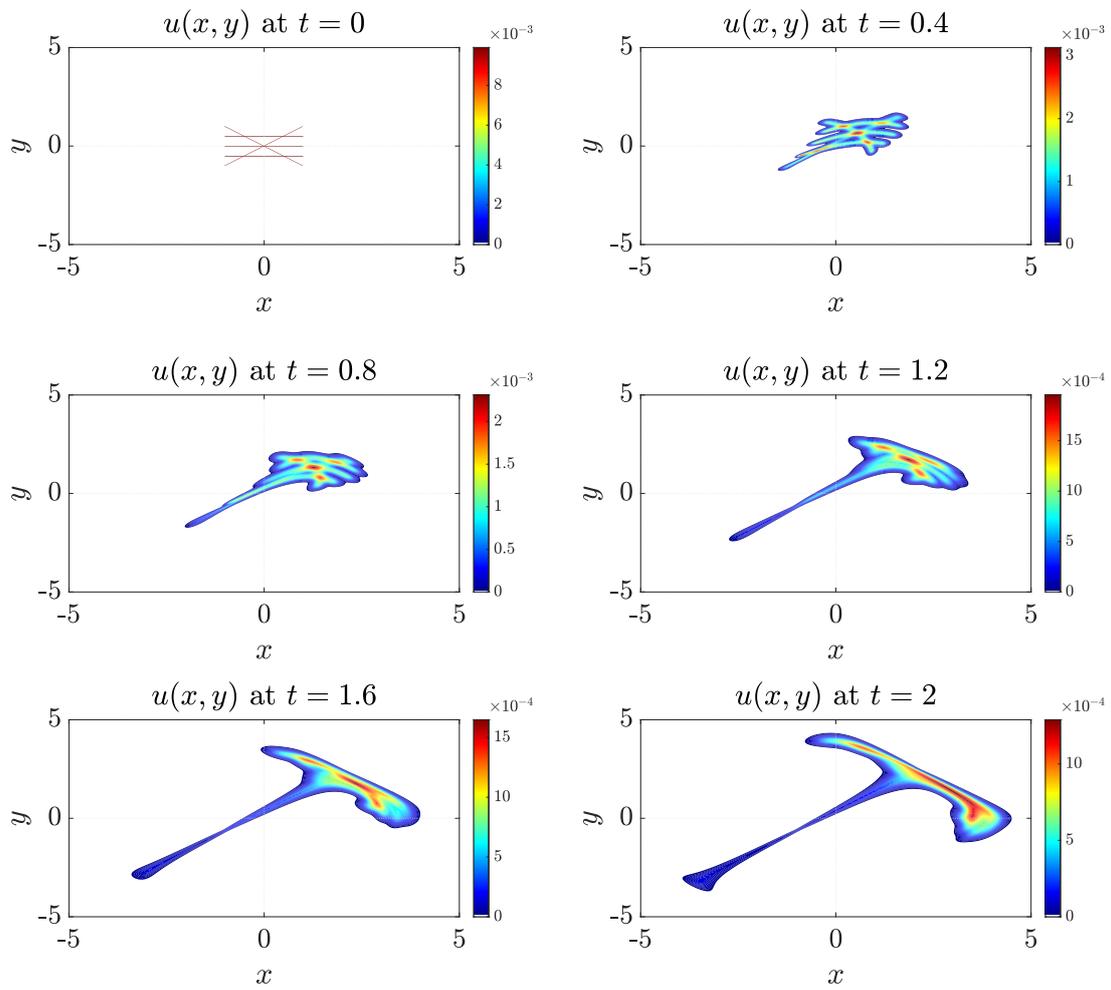

Figure 4: **Example 2** — The evolution of particles' concentration with $\beta = 10$ and smooth velocity field

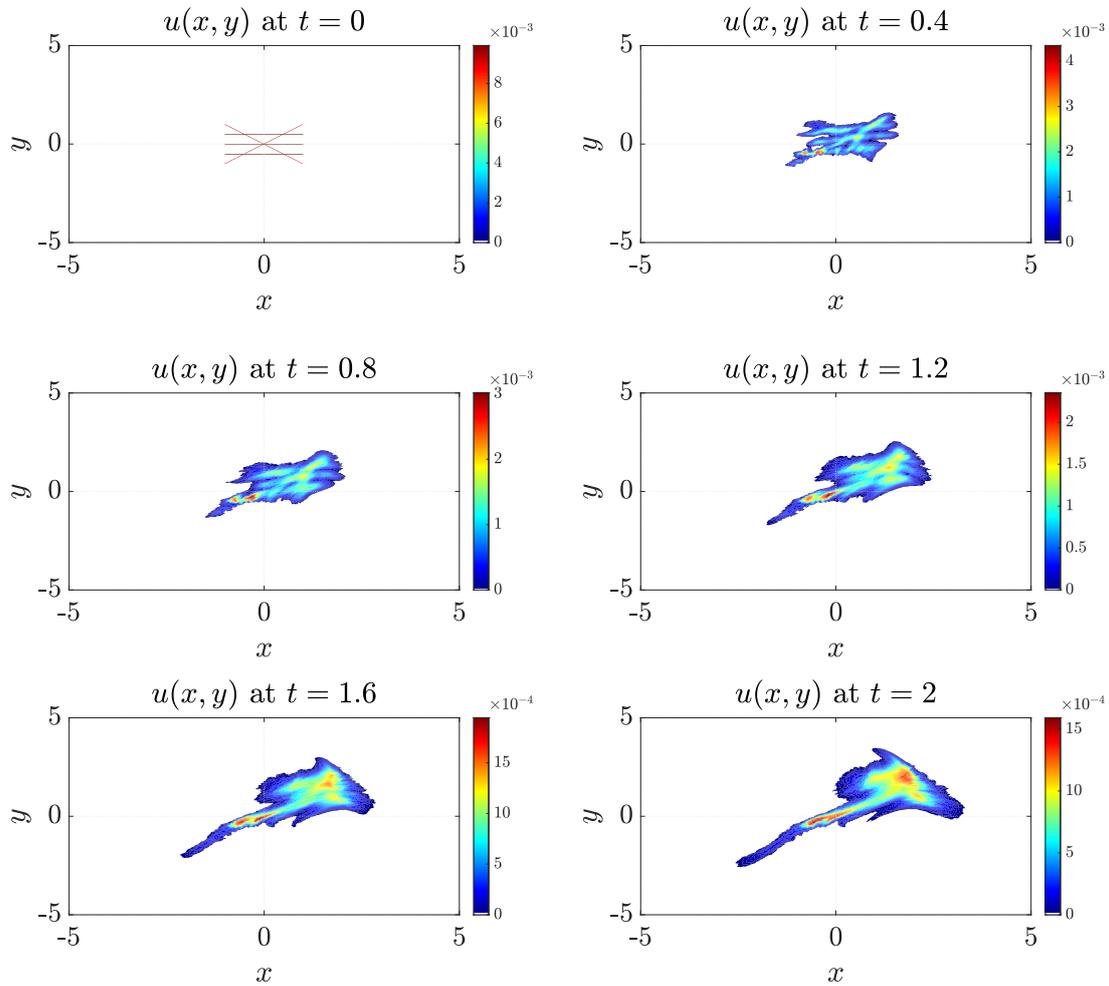

Figure 5: **Example 2** — The evolution of particles' concentration with $\beta = 10$ and noisy velocity field.

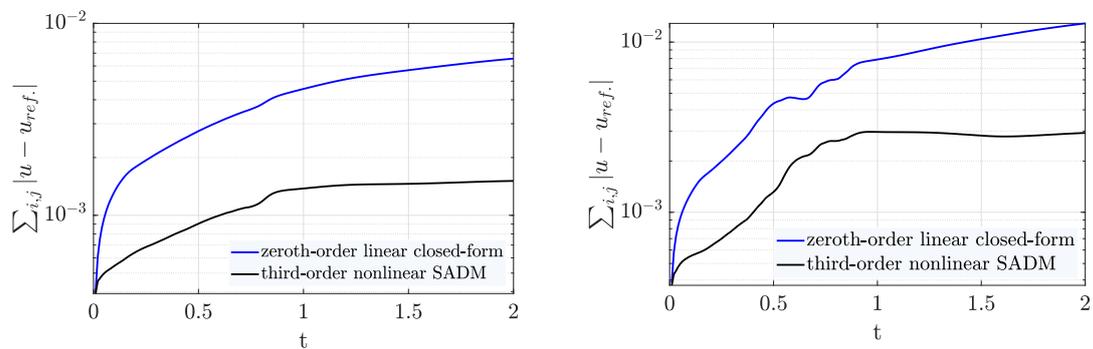

Figure 6: **Example 2** — Comparison of errors between the zeroth-order linear closed-form discretization and the third-order nonlinear SADM discretization for $\beta = 750$ (left) and $\beta = 1000$ (right).

## 4.3 Example 3: 2D Navier-Stokes Equations

The Navier-Stokes equations in the streamfunction-vorticity ($\psi$-$\omega$) formulation for 2D incompressible flow are given as follows:

$$\frac{\partial \omega}{\partial t} + \frac{\partial \psi}{\partial y}\frac{\partial \omega}{\partial x} - \frac{\partial \psi}{\partial x}\frac{\partial \omega}{\partial y} = \frac{1}{Re}\left(\frac{\partial^2 \omega}{\partial x^2} + \frac{\partial^2 \omega}{\partial y^2}\right),$$
$$\frac{\partial^2 \psi}{\partial x^2} + \frac{\partial^2 \psi}{\partial y^2} = -\omega. \qquad (35)$$

Here, $\omega$ denotes the vorticity, defined as $\omega = \frac{\partial v}{\partial x} - \frac{\partial u}{\partial y}$, and $Re$ represents the Reynolds number, expressed as $Re = \frac{UW}{\nu}$, where $U$ is the characteristic velocity, $W$ is the characteristic length (which is the channel width), and $\nu$ is the kinematic viscosity. The streamfunction $\psi$ is related to the velocity components through the relationships $u = \frac{\partial \psi}{\partial y}$, and $v = -\frac{\partial \psi}{\partial x}$ where $u$ and $v$ are the normalized velocity components in the $x$ and $y$ directions, respectively.

For the test cases considered in this study, we refer to Fig. 7 which provides a comprehensive schematic that encapsulates all the test cases. Furthermore, the analyses are conducted at $Re = 100$ with the normalized inlet velocity $u_{inlet} = 1$.

The second-order Strang splitting method was employed to solve the governing equations. The advection operator was handled using the method of characteristics, while the diffusion operator was resolved using the zeroth-order temporal-ODE scheme, as the primary focus of this example is the steady-state solution. The computational domain was discretized with 100 nodes along the $y$-direction and 400 nodes along the $x$-direction. Additionally, a fixed time step of $\Delta t = 0.01$ was utilized for all cases examined in this study.

Fig. 8 presents the $x$-direction velocity profiles, $u(x_c, y)$, at various cross-sections for all the studied cases. Specifically, the cross-sections are chosen as $x_{c,i} = i\frac{W}{4}$, $i = 1,\dots,8$, and $x_{c,i} = (i-5)W$, $i = 8,\dots,15$. Additionally, Fig. 9 illustrates the stream functions, $\psi(x, y)$, for all the configurations considered in this study.

In Fig. 10, a comparison is presented between the solutions obtained using the zeroth-order temporal-ODE scheme and the classic ADI method (detailed in **Supplementary Information, Section 7**). For each configuration, both methods are applied under identical simulation settings, including the same number of spatial nodes and iterations for solving the stream function equation (set to one iteration). The comparison focuses on the following metric, which measures the transition to a steady-state: $\sum_{i,j}\left(u^{n+1} - u^n\right)^2$.

As shown in Fig. 10, the zeroth-order temporal-ODE scheme significantly outperforms the classic ADI method in terms of the number of iterations required to achieve a steady-state. Furthermore, it is evident that the temporal-ODE scheme demonstrates superior stability and robustness for the reverse-flow configurations analyzed in this study.

Our findings indicate that the newly developed directional-ODE approach exhibits exceptional stability, even at high Reynolds numbers and larger time steps. Nevertheless, a comprehensive evaluation of its performance and accuracy under such conditions is deferred to future investigations.

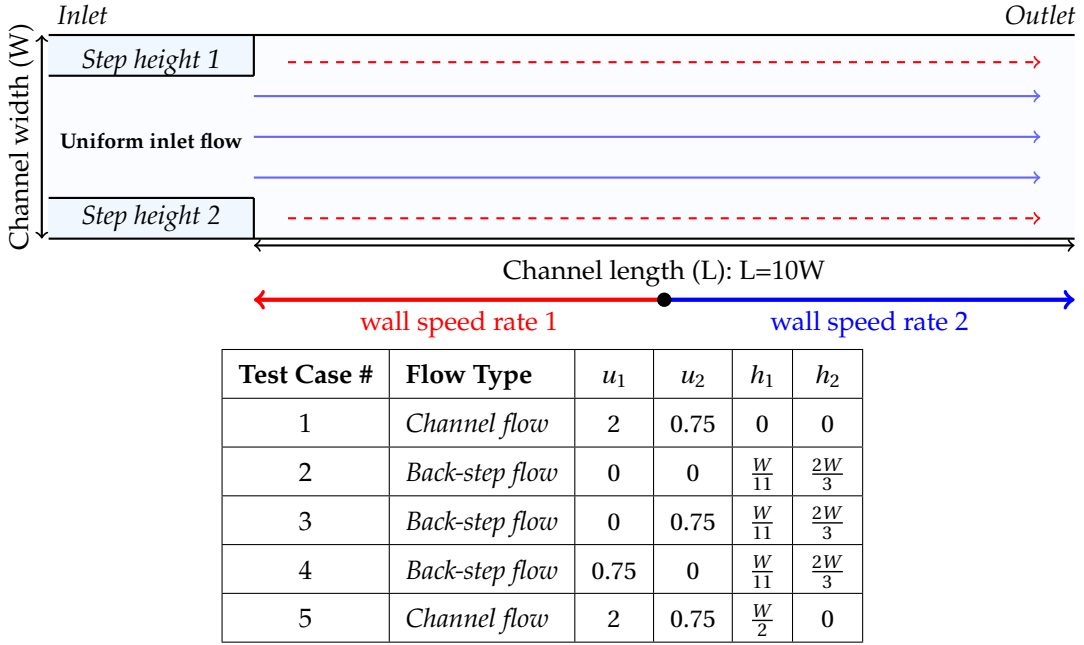

Figure 7: **Example 3** — Flow configuration schematic (top) and test case parameters (bottom): $u_1$ is **wall speed rate 1**, $u_2$ is **wall speed rate 2**, $h_1$ is **step height 1**, and $h_2$ is **step height 2**

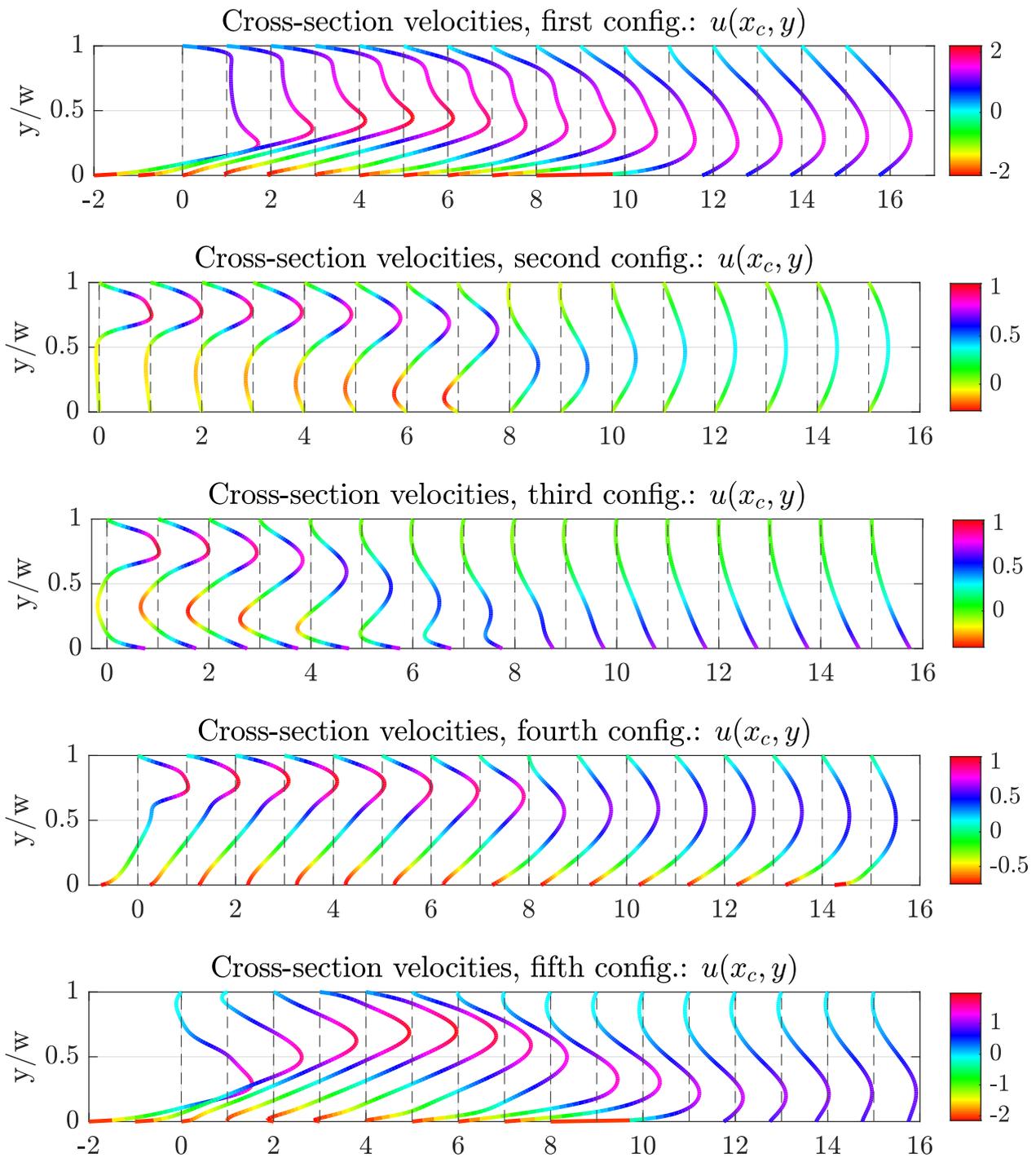

Figure 8: **Example 3** — Cross-section velocities for the considered configurations.

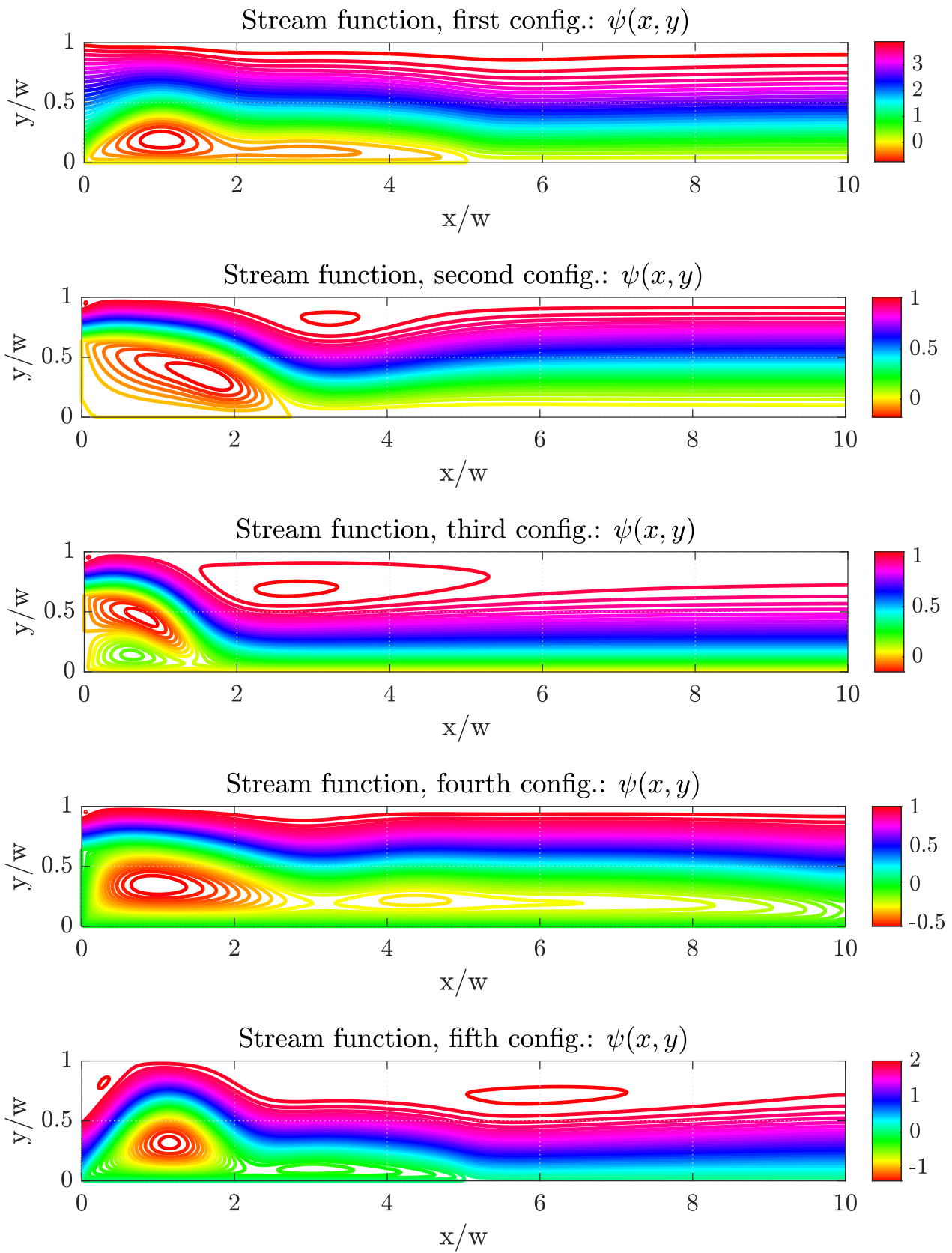

Figure 9: **Example 3** — Stream functions for the considered configurations.

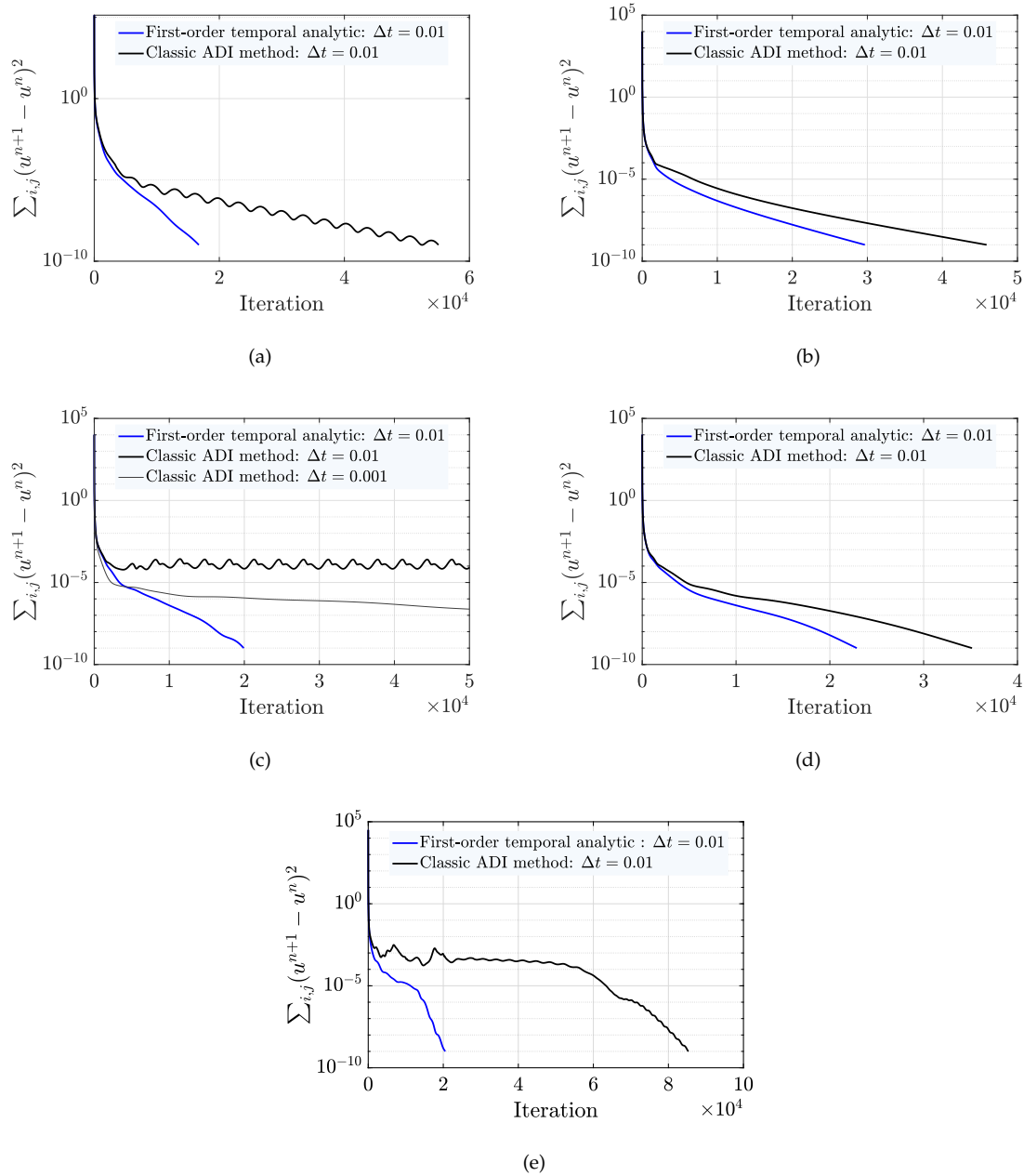

Figure 10: **Example 3**– Comparison between the zeroth-order temporal-ODE scheme and the classic ADI method– (a) to (e) represent first to fifth configurations considered respectively.

## 5 Conclusion

In this work, the directional-ODE discretization approach was introduced as a novel paradigm for discretizing PDEs, interpreting discrete schemes as directional ODEs in either the spatial or temporal domain. Focusing on ADEs, a general framework for implementing the directional-ODE discretization approach was presented, with considerations for nonlinear advection and diffusion. Sufficient stability conditions were derived, particularly for cases involving nonlinear representative ODEs. Moreover, a comprehensive stability analysis of the variable-order temporal-ODE scheme demonstrated the unconditional stability of the approach. Scenarios involving uncertain diffusion coefficients were also briefly explored, underscoring the versatility of directional-ODE approach as a powerful tool for realistic stochastic analyses.

The primary aim of this study was to introduce the general concept of the directional-ODE discretization approach, which justifies the concise and focused scope of this work. By highlighting its generality and versatility, this study establishes a strong foundation for future research and opens numerous avenues for further development.

Potential extensions include exploring higher-order spatial discretization for the diffusion operator while dealing with temporal-ODE scheme, developing high-order analytic discrete formulas for ADEs involving multiple uncertain parameters with various associated probability distributions, implementing adaptive time-step schemes, deriving high-order solutions for the advection operator, and thorough examination of predictor-corrector algorithms for temporal/spatial-ODE schemes applied directly to the ADE without splitting techniques (detailed in **Supplementary Information, Section 6**).

In addition, the central concept of the directional-ODE discretization approach holds significant promise for applications to other classes of PDEs and the analysis of ADEs without employing operator-splitting techniques. Our findings suggest that for such cases, a fully nonlinear representative ODE can be derived and solved using semi-analytic methods, such as SADM. However, a rigorous stability analysis of these representative ODEs and a thorough comparison with the high-order splitting method presented in this study (or with available analytical solutions) will be essential to validate the accuracy and practicality of such alternative approaches.

## Data and Code Availability

This work is theoretical, and as such, no experimental data is provided. However, the computational codes are publicly available in link. References to these codes are included throughout the supplementary information.

[24] C. Zoppou and J. H. Knight. Analytical solution of a spatially variable coefficient advection–diffusion equation in up to three dimensions. *Applied Mathematical Modelling*, 23(9):667–685, 1999.

# Supplementary Information: A Directional-ODE Framework for Discretization of Advection-Diffusion Equations

**Amin Jafarimoghaddam**⋆[5], **Manuel Soler**† and **Irene Ortiz**†

⋆,†*Department of Aerospace Engineering, Universidad Carlos III de Madrid. Avenida de la Universidad, 30, Leganes, 28911 Madrid, Spain* [6]

# 6 The Advection-Diffusion Equation: A General Nonlinear Model

Let us consider the following general nonlinear advection-diffusion equation (ADE):

$$\frac{\partial u(\mathbf{x},t)}{\partial t} + \mathbf{v}(\mathbf{x},t,u) \cdot \nabla u(\mathbf{x},t) = \nabla \cdot \left( \bar{d}(\mathbf{x},t,u) \nabla u(\mathbf{x},t) \right) + f(\mathbf{x},t,u). \tag{36}$$

where $\mathbf{x}$ denotes the coordinate system in a vector notation, i.e., $\mathbf{x} := (x,y,z)$ in a conventional Cartesian framework, $u$ represents the scalar quantity of interest (such as concentration), $\mathbf{v} := (\mathbf{v}_x, \mathbf{v}_y, \mathbf{v}_z)$ is the velocity field associated with advection of $u$, $\bar{d}$ is the scalar isotropic diffusion coefficient, and $f$ is the system input/source term which can be attributed to various scenarios[7].

Expanding the diffusion term, the above ADE becomes:

$$\frac{\partial u}{\partial t} + \mathbf{v} \cdot \nabla u = \frac{\partial \bar{d}}{\partial u} |\nabla u|^2 + \nabla \bar{d} \cdot \nabla u + \bar{d} \Delta u + f(\mathbf{x},t,u). \tag{37}$$

where $\Delta := \nabla \cdot (\nabla u)$ is the Laplacian of $u$, and $|\nabla u|^2 := \nabla u \cdot \nabla u$.

Specifically, the nonlinear diffusion term breaks into the following terms sub-terms:

$$\nabla \cdot \left( \bar{d}(\mathbf{x},t,u) \nabla u(\mathbf{x},t) \right) = \overbrace{\frac{\partial \bar{d}}{\partial u} |\nabla u|^2}^{\text{induced nonlinear advection}} + \overbrace{\nabla \bar{d} \cdot \nabla u}^{\text{induced linear advection}} + \overbrace{\bar{d} \Delta u}^{\text{pure diffusion}}. \tag{38}$$

In other words, the nonlinear diffusion term can be decomposed into a standard diffusion term with a direct Laplace operator, along with two additional terms involving linear and nonlinear gradient operators that account for linearly and nonlinearly shifted diffusion effects. From a computational perspective, these terms can be interpreted as additional advection components in the system. In the same section, we present a detailed comparison of two solution approaches: (1) a fully implicit method applied directly to the nonlinear diffusion term, and (2) an implicit splitting approach, where the induced linear and nonlinear advection terms are incorporated into an advection operator and solved using an operator splitting technique.

---

[5]⋆Corresponding Author.

[6]Email addresses: ajafarim@pa.uc3m.es (A. Jafarimoghaddam), masolera@ing.uc3m.es (M. Soler), irortiza@ing.uc3m.es (I. Ortiz)

[7]This document outlines a general 3D framework; however, the provided codes are limited to 1D and 2D for simplicity and user-friendliness. Our near-term goal is to develop an automated library for solving n-dimensional ADE systems, accommodating coupled equations with general initial and boundary conditions.

Reformatting the ADE, we write:

$$\frac{\partial u}{\partial t} + \mathbf{V} \cdot \nabla u = \bar{d} \Delta u + f(\mathbf{x}, t, u). \tag{39}$$

where, $\mathbf{V} := \mathbf{v} - \frac{\partial \bar{d}}{\partial u} \nabla u - \nabla \bar{d}$.

In a 3D Cartesian framework, the ADE is expanded as:

$$\frac{\partial u}{\partial t} + \mathbf{V}_x \frac{\partial u}{\partial x} + \mathbf{V}_y \frac{\partial u}{\partial y} + \mathbf{V}_z \frac{\partial u}{\partial z} = \bar{d}\left(\frac{\partial^2 u}{\partial x^2} + \frac{\partial^2 u}{\partial y^2} + \frac{\partial^2 u}{\partial z^2}\right) + f(x, y, z, t, u). \tag{40}$$

$$\begin{aligned}
\mathbf{V}_x &= \mathbf{v}_x(x, y, z, t, u) - \frac{\partial \bar{d}}{\partial u}\frac{\partial u}{\partial x} - \frac{\partial \bar{d}}{\partial x}, \\
\mathbf{V}_y &= \mathbf{v}_y(x, y, z, t, u) - \frac{\partial \bar{d}}{\partial u}\frac{\partial u}{\partial y} - \frac{\partial \bar{d}}{\partial y}, \\
\mathbf{V}_z &= \mathbf{v}_z(x, y, z, t, u) - \frac{\partial \bar{d}}{\partial u}\frac{\partial u}{\partial z} - \frac{\partial \bar{d}}{\partial z}.
\end{aligned} \tag{41}$$

For simplicity, we consider the following Dirichlet boundary conditions:

$$\begin{aligned}
u(0, y, z, t) &= u_{0,x}, \quad u(L_x, y, z, t) = u_{L_x, x}, \\
u(x, 0, z, t) &= u_{0,y}, \quad u(0, L_y, z, t) = u_{L_y, y}, \\
u(x, y, 0, t) &= u_{0,z}, \quad u(x, y, L_z, t) = u_{L_z, z}.
\end{aligned} \tag{42}$$

Moreover, the initial condition is defined as $u(x, y, z, 0) =: u^0$.

## 6.1 Operator Splitting Methods for Advection-Diffusion Equations

Double-step advection-diffusion splitting algorithm is considered to solve the ADE on $[t_n, t_{n+1}]$:

$$u^{n+1} = \mathscr{A}_{\frac{\Delta t}{2}} \mathscr{D}_{\Delta t} \mathscr{A}_{\frac{\Delta t}{2}} u^n \tag{43}$$

where $\mathscr{A}$ and $\mathscr{D}$ are the advection and diffusion operators, respectively.

## 6.2 Nonlinear Diffusion: Comparison Between Two Different Solution Approaches

Let us consider the following 1D heat equation with a nonlinear diffusion coefficient:

$$\frac{\partial u}{\partial t} = \frac{\partial}{\partial x}\left(\bar{d}(u)\frac{\partial u}{\partial x}\right). \tag{44}$$

The above equation can be re-written as:

$$\frac{\partial u}{\partial t} = \left(\frac{\partial \bar{d}}{\partial u}\right)\left(\frac{\partial u}{\partial x}\right)^2 + \bar{d}\frac{\partial^2 u}{\partial x^2}. \tag{45}$$

Here, we solve the above equation employing the following approaches:

**1)** *Implicit discretization without operator splitting technique*: In this case, the nonlinear term $\left(\frac{\partial u}{\partial x}\right)^2$ is initially linearized as $\left(\frac{\partial u}{\partial x}\right)^2 \approx 2\left(\frac{\partial u}{\partial x}\right)|^n \left(\frac{\partial u}{\partial x}\right)|^{n+1} - \left(\frac{\partial u}{\partial x}|^n\right)^2$. Therefore, we write the full implicit scheme (using the central difference scheme) as:

$$\frac{u_i^{n+1} - u_i^n}{\Delta t} = \frac{\partial \bar{d}}{\partial u}\Big|_i^n\left(\frac{\partial u}{\partial x}\Big|^n\frac{\partial u}{\partial x}\Big|^{n+1} - \left(\frac{\partial u}{\partial x}\Big|^n\right)^2\right) + \bar{d}(u_i^n)\frac{u_{i-1}^{n+1} - 2u_i^{n+1} + u_{i+1}^{n+1}}{(\Delta x)^2}. \tag{46}$$

where $\frac{\partial u}{\partial x}|^{n+1} = \frac{u_{i+1}^{n+1} - u_{i-1}^{n+1}}{2\Delta x}$.

The above implicit discretization scheme can be rewritten as:

$$au_{i-1}^{n+1} + bu_i^{n+1} + cu_{i+1}^{n+1} = d. \tag{47}$$

where,

$$\begin{aligned}
a &= -\frac{\Delta t}{\Delta x}\frac{\partial u}{\partial x}\Big|^n \frac{\partial \bar{d}}{\partial u}\Big|_i^n + \lambda, \\
b &= -1 - 2\lambda, \\
c &= \frac{\Delta t}{\Delta x}\frac{\partial u}{\partial x}\Big|^n \frac{\partial \bar{d}}{\partial u}\Big|_i^n + \lambda, \\
d &= -u_i^n + \left(\frac{\partial u}{\partial x}\Big|^n\right)^2 \frac{\partial \bar{d}}{\partial u}\Big|_i^n \Delta t.
\end{aligned} \tag{48}$$

In above, $\lambda := \Delta t \frac{\bar{d}(u_i^n)}{(\Delta x)^2}$ is the CFL factor, $\frac{\partial u}{\partial x}|^n = \frac{u_{i+1}^n - u_{i-1}^n}{2\Delta x}$. Moreover, $\bar{d}(u)$ and $\frac{\partial \bar{d}}{\partial u}$ are computed at $t^n$.

The above triangular system is solved using TDMA.

**2)** *Implicit discretization with operator splitting technique*. In this case, we initially write:

$$\frac{\partial u}{\partial t} - \mathbf{V}_x\frac{\partial u}{\partial x} = \bar{d}\frac{\partial^2 u}{\partial x^2}, \quad \mathbf{V}_x := \left(\frac{\partial \bar{d}}{\partial u}\right)\left(\frac{\partial u}{\partial x}\right). \tag{49}$$

Next, we apply a single-step advection-diffusion splitting algorithm to solve the advection-diffusion equation (ADE) over the interval $[t_n, t_{n+1}]$:

$$u^{n+1} = \mathscr{D}_{\Delta t}\mathscr{A}_{\Delta t} u^n, \tag{50}$$

where the operators $\mathscr{A}_{\Delta t}$ and $\mathscr{D}_{\Delta t}$ are defined as:

**a. Advection Step (Method of Characteristics):**

$$\mathscr{A}_{\Delta t} u^n : \frac{\partial u}{\partial t} + \mathbf{V}_x \frac{\partial u}{\partial x} = 0, \quad \text{with initial condition } u^n \to u^*. \tag{51}$$

**b. Diffusion Step (Implicit Discretization[8]):**

$$\mathscr{D}_{\Delta t} u^* : \frac{\partial u}{\partial t} = \bar{d}\frac{\partial^2 u}{\partial x^2}, \quad \text{with initial condition } u^* \to u^{n+1}. \tag{52}$$

The code available at link solves a 1D heat equation with user-defined arbitrary nonlinear diffusion coefficients, applying the solution strategies discussed above and comparing their performance. The results demonstrate a perfect match between the two methods. Notably, the *implicit discretization with the operator splitting technique* is shown to be unconditionally stable. In contrast, for strongly nonlinear diffusion coefficients, the *implicit discretization scheme without operator splitting technique* proves ineffective in dealing with the term $\left(\frac{\partial \bar{d}}{\partial u}\right)\left(\frac{\partial u}{\partial x}\right)^2$.

# 7 Elucidation of the Solution Methodology in a full 3D Framework

In this section, we document the solution methodology used for a full 3D framework. However, since the examples in the present paper are at most 2D, we also provide open-source 2D codes with detailed instructions on their usage.

## 7.1 Computational Domain

Let us initially discretize the computational domain as:

$$\begin{aligned}
& x \in [0, Lx], \quad y \in [0, L_y], \quad z \in [0, Lz], \quad t \in [0, T], \\
& \Delta x := \frac{L_x}{N_x - 1} = x_{i+1} - x_i, \quad \Delta y := \frac{L_y}{N_y - 1} = y_{j+1} - y_j, \quad \Delta z := \frac{L_z}{N_z - 1} = z_{k+1} - z_k, \\
& i = 0, ..., N_x - 1, \quad j = 0, ..., N_y - 1, \quad k = 0, ..., N_z - 1, \\
& \Delta t := \frac{T}{N_t - 1} = t_{n+1} - t_n, \quad n = 0, ..., M - 1.
\end{aligned} \tag{53}$$

---

[8]The implicit discretization for the pure diffusion equation is given by $\frac{u_i^{n+1} - u_i^n}{\Delta t} = \bar{d}(u_i^n)\frac{u_{i-1}^{n+1} - 2u_i^{n+1} + u_{i+1}^{n+1}}{(\Delta x)^2}$, resulting in a triangular system that can be efficiently solved using the TDMA.

In above, $x_0 := 0, x_{N_x} := L_x, y_0 := 0, y_{N_y} := L_y, z_0 := 0, z_{N_z} := L_z$, and $t_0 := 0, t_M := T$.

## 7.2 Advection Operator: The Method of Characteristics

For the advection operator, we exploit the method of characteristics to solve the following advection equation on $[t_n, t_{n+\frac{1}{2}}]$:

$$\frac{\partial u}{\partial t} + \mathbf{V}_x \frac{\partial u}{\partial x} + \mathbf{V}_y \frac{\partial u}{\partial y} + \mathbf{V}_z \frac{\partial u}{\partial z} = 0. \tag{54}$$

The method of characteristics involves solving the following system ODEs to trace the characteristics:

$$\frac{dx}{dt} = \mathbf{V}_x, \quad x(t_n) = x_0, \quad \frac{dy}{dt} = \mathbf{V}_y, \quad y(t_n) = y_0, \quad \frac{dz}{dt} = \mathbf{V}_z, \quad z(t_n) = z_0. \tag{55}$$

where $(x, y, z) =: \mathbf{x}$ describes the path of the particles along the characteristics, and the solution is constant along these paths:

$$u(x, y, z, t) = u(x_0, y_0, z_0, t_n). \tag{56}$$

where $u(x_0, y_0, z_0, t_n)$ is the initial condition.

For a general semi-analytic solution, one can apply methods such as the SADM. On using the SADM over $[t_n, t_{n+\frac{1}{2}}]$, the characteristic solution is obtained as:

$$\begin{aligned}
x_0 &= x - \mathscr{L}^{-1}[\mathbf{V}_x] := x - \int_{t_n}^{t_{n+\frac{1}{2}}} \mathbf{A}_{0,x} dt - \int_{t_n}^{t_{n+\frac{1}{2}}} \mathbf{A}_{1,x} dt - ..., \\
y_0 &= y - \mathscr{L}^{-1}[\mathbf{V}_y] := y - \int_{t_n}^{t_{n+\frac{1}{2}}} \mathbf{A}_{0,y} dt - \int_{t_n}^{t_{n+\frac{1}{2}}} \mathbf{A}_{1,y} dt - ..., \\
z_0 &= z - \mathscr{L}^{-1}[\mathbf{V}_z] := z - \int_{t_n}^{t_{n+\frac{1}{2}}} \mathbf{A}_{0,z} dt - \int_{t_n}^{t_{n+\frac{1}{2}}} \mathbf{A}_{1,z} dt - ....
\end{aligned} \tag{57}$$

where $\mathbf{A}_{j,\mathbf{x}}$ represents Adomian terms in $x, y$, and $z$ directions.

Notably, depending on the desired accuracy, one may also adopt the Euler approximation as the simplest approach to solving the above ODE on the interval $[t_n, t_{n+\frac{1}{2}}]$, resulting in: $x_0 = x - \mathbf{V}_x \frac{\Delta t}{2}$, $y_0 = y - \mathbf{V}_y \frac{\Delta t}{2}$, and $z_0 = z - \mathbf{V}_z \frac{\Delta t}{2}$.

Therefore, the general solution to the advection operator involves the following coordinate mapping at each time step:

$$\mathscr{A}_{\frac{\Delta t}{2}} u^n := u\left(x - \mathscr{L}^{-1}[\mathbf{V}_x], y - \mathscr{L}^{-1}[\mathbf{V}_y], z - \mathscr{L}^{-1}[\mathbf{V}_z]\right). \tag{58}$$

The coordinate mapping described above can be achieved using various interpolation methods (such as linear interpolation), which are readily available in programming languages such as MATLAB.

## 7.3 Diffusion Operator: Temporal-ODE Discretization Scheme

The diffusion operator involves solving the following heat equation:

$$\frac{\partial u}{\partial t} = \bar{d}\left(\frac{\partial^2 u}{\partial x^2} + \frac{\partial^2 u}{\partial y^2} + \frac{\partial^2 u}{\partial z^2}\right) + f(x,y,z,t,u). \tag{59}$$

To develop a temporal-ODE scheme of the diffusion operator, we initially write the heat equation as:

$$\begin{aligned}\frac{du_{i,j,k}}{dt} &= \bar{d}\left(\frac{u_{i-1,j,k}(t) + u_{i+1,j,k}(t)}{(\Delta x)^2} + \frac{u_{i,j-1,k}(t) + u_{i,j+1,k}(t)}{(\Delta y)^2} + \frac{u_{i,j,k-1}(t) + u_{i,j,k+1}(t)}{(\Delta z)^2}\right) \\ &\quad - 2u_{i,j,k}(t)\left(\frac{\bar{d}}{(\Delta x)^2} + \frac{\bar{d}}{(\Delta y)^2} + \frac{\bar{d}}{(\Delta z)^2}\right) + f(x_{i,j,k}, y_{i,j,k}, z_{i,j,k}, t_n, u_{i,j,k}^n).\end{aligned} \tag{60}$$

Defining:

$$\bar{a} := \bar{d}\left(\frac{1}{(\Delta x)^2} + \frac{1}{(\Delta y)^2} + \frac{1}{(\Delta z)^2}\right),$$

$$\mathcal{U}(\tau) := \sum_{p=0}^{P} a_p \tau^p := \frac{u_{i-1,j,k}(\tau) + u_{i+1,j,k}(\tau)}{(\Delta x)^2} + \frac{u_{i,j-1,k}(\tau) + u_{i,j+1,k}(\tau)}{(\Delta y)^2} + \frac{u_{i,j,k-1}(\tau) + u_{i,j,k+1}(\tau)}{(\Delta z)^2}, \tag{61}$$

$\tau := t - t_n, \forall t \in [t_n, t_{n+1}] \to \tau \in [0, \Delta t]$.

The ODE for $u_{i,j,k}(t)$ can be written as:

$$\frac{du}{d\tau} = \bar{d}\left(\sum_{p=0}^{P} a_p \tau^p\right) - 2\bar{a}u + s. \tag{62}$$

The solution reads:

$$u_{i,j,k}(\tau) = ce^{-2\bar{a}\tau} + \frac{\bar{d}}{\bar{a}} \sum_{p=0}^{P} \frac{a_p}{2} \sum_{q=0}^{p} \frac{\tau^{p-q}(-1)^q p!}{(2\bar{a})^q (p-q)!} + \frac{s}{2\bar{a}}. \tag{63}$$

where, $s := f(x_{i,j,k}, y_{i,j,k}, z_{i,j,k}, t_n, u_{i,j,k}^n)$ is known at $t_n$, and:

$$c = u_{i,j,k}^n - \frac{s}{2\bar{a}} - \frac{\bar{d}}{\bar{a}} \sum_{p=0}^{P} \frac{a_p}{2} \frac{(-1)^p p!}{(2\bar{a})^p}. \tag{64}$$

Therefore, we can write the closed-form solution as:

$$u_{i,j,k}(\tau) = \left(u_{i,j,k}^n - \frac{s}{2\bar{a}} - \frac{\bar{d}}{\bar{a}} \sum_{p=0}^{P} \frac{a_p}{2} \frac{(-1)^p p!}{(2\bar{a})^p}\right) e^{-2\bar{a}\tau} + \frac{\bar{d}}{\bar{a}} \sum_{p=0}^{P} \frac{a_p}{2} \sum_{q=0}^{p} \frac{\tau^{p-q}(-1)^q p!}{(2\bar{a})^q (p-q)!} + \frac{s}{2\bar{a}}. \tag{65}$$

In above, $a_p$ are calculated according to an intermediate sampling over $[0, \Delta t]$.

## 7.4 The Predictor-Corrector Algorithm

**Predictor Step:** We initially predict intermediate values of $\mathcal{U}(\tau)$ using the update formula associated with $P = 0$:

$$u_{i,j,k}^{predictor}(\tau) =: \mathcal{F}_{predictor}(u_{i,j,k}^n, \tau) = \left(u_{i,j,k}^n - \frac{s}{2\bar{a}} - \frac{\bar{d}}{\bar{a}} \frac{a_0}{2}\right) e^{-2\bar{a}\tau} + \frac{\bar{d}}{\bar{a}} \frac{a_0}{2} + \frac{s}{2\bar{a}}. \tag{66}$$

Therefore, the sample points are:

$$\mathscr{S} := \{(0, \mathscr{U}^n_{predictor}), (T_1, \mathscr{U}^{n+\frac{1}{N}}_{predictor}), \ldots, (T_{N-1}, \mathscr{U}^{n+\frac{(N-1)}{N}}_{predictor}), (T_N, \mathscr{U}^{n+1}_{predictor})\}. \tag{67}$$

Here, $N \in \mathbb{N}$ denotes the order of the required polynomial approximation for the neighboring nodes. Moreover, $T_k$, for $k = 1, \ldots, N$, represents the intermediate time steps at which the samples are collected. For example, uniform time stepping is given by:

$$T_K = K \frac{\Delta t}{N}, \quad K = 1, \ldots, N. \tag{68}$$

Notably, as discussed in the paper, for increased accuracy, non-uniform sampling points, such as Chebyshev nodes, can be employed.

**Corrector Step:** Using $N+1$ sample points, we can symbolically solve the following system to determine the polynomial coefficients $a_p$ for $p = 0, 1, \ldots, N = P$, in terms of the predicted neighboring nodes:

$$\overbrace{\begin{pmatrix} 1 & 0 & 0 & \cdots & 0 \\ 1 & T_1^1 & T_1^2 & \cdots & T_1^P \\ 1 & T_2^1 & T_2^2 & \cdots & T_2^P \\ \vdots & \vdots & \vdots & \ddots & \vdots \\ 1 & T_{N-1}^1 & T_{N-1}^2 & \cdots & T_{N-1}^P \\ 1 & T_N^1 & T_N^2 & \cdots & T_N^P \end{pmatrix}}^{\mathbf{M}} \overbrace{\begin{pmatrix} a_0 \\ a_1 \\ a_2 \\ \vdots \\ a_{P-1} \\ a_P \end{pmatrix}}^{\mathbf{a}} = \overbrace{\begin{pmatrix} \mathscr{U}^n_{predictor} \\ \mathscr{U}^{n+\frac{1}{N}}_{predictor} \\ \mathscr{U}^{n+\frac{2}{N}}_{predictor} \\ \vdots \\ \mathscr{U}^{n+\frac{N-1}{N}}_{predictor} \\ \mathscr{U}^{n+1}_{predictor} \end{pmatrix}}^{\mathbf{U}}. \tag{69}$$

Therefore, the solution reads: $\mathbf{a} := \mathbf{M}^{-1}\mathbf{U}$.

Next, we exploit $a_p, p = 0, 1, \ldots, P$ into the corresponding $P^{th}$-order closed-form solution to compute $u_{i,j,k}^{n+1}$, update the coefficients and repeat the process for the correction step.

Below, the **predictor-corrector algorithm** is sketched for a general multi-stage $P^{th}$-order polynomial correction:

**Algorithm 1** Predictor-Corrector Algorithm: Marching from $u_{i,j,k}^n$ to $u_{i,j,k}^{n+1}$:

1: **Step 1: Predictor Stage**
2: $\mathcal{U}_{\text{predictor}}^{n+\frac{K}{N}} \leftarrow u_{i,j,k}^{\text{predictor}}(K\frac{\Delta t}{N}) \leftarrow \mathcal{F}_{\text{predictor}}(u_{i,j,k}^n, K\frac{\Delta t}{N}), \quad K = 1, 2, ..., N.$
3: **Step 2: Update Coefficients**
4: $\mathbf{a} \leftarrow (a_0, a_1, a_2, ..., a_P)^T \leftarrow \mathbf{M}^{-1}\mathbf{U}$
5: **Step 3: Corrector Stage**
6: **for** iteration = 1 to max(iteration) **do**
7: $\quad \mathcal{U}_{\text{predictor}}^{n+\frac{K}{N}} \leftarrow u_{i,j,k}^{\text{corrector}}(K\frac{\Delta t}{N}) \leftarrow \mathcal{F}_{\text{corrector}}(u_{i,j,k}^n, \mathbf{a}, K\frac{\Delta t}{N}), \quad K = 1, 2, ..., N.$
8: $\quad \mathbf{a} \leftarrow (a_0, a_1, a_2, ..., a_P)^T \leftarrow \mathbf{M}^{-1}\mathbf{U}$
9: $\quad$ **if** $\left| u_{i,j,k}^{\text{corrector}}(\Delta t)|_{\text{iteration}} - u_{i,j,k}^{\text{corrector}}(\Delta t)|_{\text{iteration-1}} \right| < \delta$ **then**
10: $\quad\quad$ **break**
11: $\quad$ **end if**
12: **end for**

## 7.5 Solution Methodology: Second-Order Example

For clarity, we explicitly present the methodology for the second-order temporal *analytic discretization* scheme.

**Predictor Step:** We initially predict intermediate values of $\mathcal{U}(\tau)$ using the update formula associated with $P = 0$:

$$u_{i,j,k}^{predictor}(\tau) =: \mathcal{F}_{predictor}(u_{i,j,k}^n, \tau) = \left(u_{i,j,k}^n - \frac{s}{2\bar{a}} - \frac{\bar{d}}{\bar{a}}\frac{a_0}{2}\right)e^{-2\bar{a}\tau} + \frac{\bar{d}}{\bar{a}}\frac{a_0}{2} + \frac{s}{2\bar{a}}. \tag{70}$$

where $\tau \in [0, \Delta t]$, and:

$$a_0 = \frac{u_{i-1,j,k}(0) + u_{i+1,j,k}(0)}{(\Delta x)^2} + \frac{u_{i,j-1,k}(0) + u_{i,j+1,k}(0)}{(\Delta y)^2} + \frac{u_{i,j,k-1}(0) + u_{i,j,k+1}(0)}{(\Delta z)^2}. \tag{71}$$

To construct a second-order temporal *analytic discretization* scheme, three sample points are required, which can be collected at equidistant intermediate time steps as follows:

$$\mathcal{S}_{\text{second-order}} := \{(0, \mathcal{U}_{predictor}^n), (\frac{\Delta t}{2}, \mathcal{U}_{predictor}^{n+\frac{1}{2}}), (\Delta t, \mathcal{U}_{predictor}^{n+1})\}. \tag{72}$$

Notably, $\mathcal{U}_{predictor}^n$ is already known at $t = t_n$ ($\tau = 0$).

Specifically, to perform the sampling process, we first compute $u_{i,j,k}^{predictor}(\frac{\Delta t}{2})$ and $u_{i,j,k}^{predictor}(\Delta t)$ for all grid points. Using these values, we then obtain $\mathcal{U}_{predictor}^{n+\frac{1}{2}} := \mathcal{U}(\frac{\Delta t}{2})$ and $\mathcal{U}_{predictor}^{n+1} := \mathcal{U}(\Delta t)$

through:

$$\mathcal{U}_{predictor} = \mathcal{U}_{i,j,k}(\tau) = \frac{u^{predictor}_{i-1,j,k}(\tau) + u^{predictor}_{i+1,j,k}(\tau)}{(\Delta x)^2} + \frac{u^{predictor}_{i,j-1,k}(\tau) + u^{predictor}_{i,j+1,k}(\tau)}{(\Delta y)^2} + \frac{u^{predictor}_{i,j,k-1}(\tau) + u^{predictor}_{i,j,k+1}(\tau)}{(\Delta z)^2}. \tag{73}$$

where $\tau = \frac{\Delta t}{2}, \Delta t$ corresponds to $\mathcal{U}^{n+\frac{1}{2}}_{predictor}, \mathcal{U}^{n+1}_{predictor}$ respectively.

**Corrector Step:** The coefficients of a second-order polynomial satisfy the following algebraic system:

$$\overbrace{\begin{pmatrix} 1 & 0 & 0 \\ 1 & \left(\frac{\Delta t}{2}\right)^1 & \left(\frac{\Delta t}{2}\right)^2 \\ 1 & (\Delta t)^1 & (\Delta t)^2 \end{pmatrix}}^{\mathbf{M}} \overbrace{\begin{pmatrix} a_0 \\ a_1 \\ a_2 \end{pmatrix}}^{\mathbf{a}} = \overbrace{\begin{pmatrix} \mathcal{U}^n_{predictor} \\ \mathcal{U}^{n+\frac{1}{2}}_{predictor} \\ \mathcal{U}^{n+1}_{predictor} \end{pmatrix}}^{\mathbf{U}}. \tag{74}$$

The symbolic solution for $a_0, a_1$, and $a_2$ reads:

$$\begin{aligned} a_0 &= \mathcal{U}^n_p, \\ a_1 &= \frac{1}{\Delta t}\left(-3\mathcal{U}^n_p + 4\mathcal{U}^{n+\frac{1}{2}}_p - \mathcal{U}^{n+1}_p\right), \\ a_2 &= \frac{1}{(\Delta t)^2}\left(2\mathcal{U}^n_p - 4\mathcal{U}^{n+\frac{1}{2}}_p + 2\mathcal{U}^{n+1}_p\right). \end{aligned} \tag{75}$$

In above, the subscript $p$ denotes the predictor values.

Next, with the above coefficients, the corrector values $u^{corrector}_{i,j,k}(\tau)$ are computed as:

$$u^{corrector}_{i,j,k}(\tau) =: \mathcal{F}_{corrector}(u^n_{i,j,k}, \mathbf{a}, \tau) = \frac{s}{2\bar{a}} + e^{-2\bar{a}\tau}\left(u^n_{i,j,k} - \frac{s}{2\bar{a}} - \frac{\bar{d}}{\bar{a}}\left(\frac{a_0}{2} - \frac{a_1}{4\bar{a}} + \frac{a_2}{4\bar{a}^2}\right)\right) + \frac{\bar{d}}{\bar{a}}\left(\frac{a_0}{2} + \frac{a_1}{2}\left(\tau - \frac{1}{2\bar{a}}\right) + \frac{a_2}{2}\left(\frac{1}{2\bar{a}^2} - \frac{\tau}{\bar{a}} + \tau^2\right)\right). \tag{76}$$

The code available at link solves general 2D ADEs with user-defined initial/boundary conditions using the second-order temporal-ODE scheme.

# 8 Proof of the Convergence for the Predictor-Corrector Algorithm

Following the discussion in **Section 3 of the Main Paper**, the stability analysis focuses on the diffusion operator ($s = 0$). To ensure the numerical scheme remains stable as $\Delta t \to \infty$, we analyze its asymptotic behavior. Specifically, we show that the discrete update formulas for an arbitrary order $P$ reduce to either the well-known Explicit Gauss-Seidel iteration (for $P = 0$) or the Fully-Implicit Iteration (for $P \neq 0$) when solving the steady-state diffusion equation:

$$\lim_{\Delta t \to \infty} u_{i,j,k}(\Delta t, P) = \text{Explicit Gauss-Seidel Update or Fully-Implicit Update}. \tag{77}$$

The Explicit Gauss-Seidel update at each grid point for $\nabla^2 u = 0$ is expressed as:

$$\text{Explicit: } u^{r+1} = \frac{1}{W} \sum_{\text{neighbors}} W_i u_i^r. \tag{78}$$

On the other hand, the Fully-Implicit Update at each grid point for $\nabla^2 u = 0$ is expressed as:

$$\text{Implicit: } u^{r+1} = \frac{1}{W} \sum_{\text{neighbors}} W_i u_i^{r+1}. \tag{79}$$

Here, the LHS contains the new update at the current iteration $r + 1$, while the RHS is a weighted sum of neighboring values, all taken from the previous iteration $r$. Moreover, $W$ and $W_i$ are weights determined by the grid spacing. It is well known that the Fully-Implicit Update is unconditionally stable. Similarly, the Gauss-Seidel method is known to be unconditionally stable for diffusion-type problems, ensuring that the error $e^r = u^r - u^*$ satisfies $\|e^{r+1}\| \leq \rho \|e^r\|$, $0 \leq \rho < 1$.

Using the notations introduced in this document, the above update formulas are equivalent to:

$$\text{Explicit Gauss-Seidel Update}: u_{i,j,k}^{n+1} = \frac{\bar{d}}{\bar{a}} \frac{\mathcal{U}^n}{2}, \quad \text{Fully-Implicit Update}: u_{i,j,k}^{n+1} = \frac{\bar{d}}{\bar{a}} \frac{\mathcal{U}^{n+1}}{2}. \tag{80}$$

where $\mathcal{U}^n$, and $\mathcal{U}^{n+1}$ contain the neighboring grid points (see Eq. (61)).

Therefore, it suffices to show that:

$$\lim_{\Delta t \to \infty} u(\Delta t, P) = \frac{\bar{d}}{\bar{a}} \cdot \frac{\mathcal{U}^{n+\delta}}{2}, \quad \delta \in \{0, 1\}. \tag{81}$$

## 8.1 Proof of the Proposition: $\lim_{\Delta t \to \infty} u(\Delta t, P) = \frac{\bar{d}}{\bar{a}} \cdot \frac{\mathcal{U}^{n+\delta}}{2}, \quad \delta \in \{0, 1\}$

For simplicity, we consider equidistant sampling, defined as $T_K = K\frac{\Delta t}{N}$, $K = 1, \ldots, N$. It can be shown that:

$$\mathbf{M}^{-1} = \text{diag}\left((\Delta t)^0, (\Delta t)^{-1}, \ldots, (\Delta t)^{-P}\right) \bar{\mathbf{M}}. \tag{82}$$

where $\bar{\mathbf{M}}$ is a coefficient matrix of size $P \times P$ that is independent of $\Delta t$.

The polynomial coefficients are computed symbolically by solving $\mathbf{a} = \mathbf{M}^{-1}\mathbf{U}$, which is equivalent to $\mathbf{a} = \text{diag}\left((\Delta t)^0, (\Delta t)^{-1}, \ldots, (\Delta t)^{-P}\right) \bar{\mathbf{M}} \mathbf{U}$.

Computing for the first few terms, we obtain:

$$P = 0: \bar{\mathbf{M}} = \begin{pmatrix} 1 \end{pmatrix}, \quad P = 1: \bar{\mathbf{M}} = \begin{pmatrix} 1 & 0 \\ -1 & 1 \end{pmatrix}, \quad P = 2: \bar{\mathbf{M}} = \begin{pmatrix} 1 & 0 & 0 \\ -3 & 4 & -1 \\ 2 & -4 & 2 \end{pmatrix}. \tag{83}$$

and so on.

Therefore, we obtain the coefficients as:

$$\begin{aligned} P = 0: \; & a_0 = \mathscr{U}^n, \\ P = 1: \; & a_0 = \mathscr{U}^n, \quad a_1 = \frac{1}{\Delta t}(\mathscr{U}^{n+1} - \mathscr{U}^n), \\ P = 2: \; & a_0 = \mathscr{U}^n, \quad a_1 = \frac{1}{\Delta t}(-3\mathscr{U}^n + 4\mathscr{U}^{n+\frac{1}{2}} - \mathscr{U}^{n+1}), \quad a_2 = \frac{1}{(\Delta t)^2}(2\mathscr{U}^n - 4\mathscr{U}^{n+\frac{1}{2}} + 2\mathscr{U}^{n+1}). \end{aligned} \tag{84}$$

Using the above coefficients, update formulas for $P = 0, 1, 2$ are obtained as:

**P = 0:**

$$\lim_{\Delta t \to \infty} u_{i,j,k}(\Delta t, P = 0) = \lim_{\Delta t \to \infty} \left( \left( u_{i,j,k}^n - \frac{\bar{d}}{\bar{a}} \frac{a_0}{2} \right) e^{-2\bar{a}\Delta t} + \frac{\bar{d}}{\bar{a}} \frac{a_0}{2} \right) = \frac{\bar{d}}{\bar{a}} \frac{a_0}{2} = \frac{\bar{d}}{\bar{a}} \frac{\mathscr{U}^n}{2}. \tag{85}$$

**P = 1:**

$$\begin{aligned} & \lim_{\Delta t \to \infty} u_{i,j,k}(\Delta t, P = 1) = \\ & \lim_{\Delta t \to \infty} \left( e^{-2\bar{a}\Delta t} \left( u_{i,j,k}^n - \frac{\bar{d}}{\bar{a}} \left( \frac{a_0}{2} - \frac{a_1}{4\bar{a}} \right) \right) + \frac{\bar{d}}{\bar{a}} \left( \frac{a_0}{2} + \frac{a_1}{2} \left( \Delta t - \frac{1}{2\bar{a}} \right) \right) \right) = \\ & \lim_{\Delta t \to \infty} \left( \frac{\bar{d}}{\bar{a}} \left( \frac{a_0}{2} + \frac{a_1}{2} \left( \Delta t - \frac{1}{2\bar{a}} \right) \right) \right) = \\ & \lim_{\Delta t \to \infty} \left( \frac{\bar{d}}{\bar{a}} \left( \frac{\mathscr{U}^n}{2} + \frac{\mathscr{U}^{n+1} - \mathscr{U}^n}{2\Delta t} \left( \Delta t - \frac{1}{2\bar{a}} \right) \right) \right) = \\ & \lim_{\Delta t \to \infty} \left( \frac{\bar{d}}{\bar{a}} \left( \frac{\mathscr{U}^n}{2} + \frac{\mathscr{U}^{n+1} - \mathscr{U}^n}{2} \right) \right) = \frac{\bar{d}}{\bar{a}} \left( \frac{\mathscr{U}^{n+1}}{2} \right). \end{aligned} \tag{86}$$

**P = 2:**

$$\begin{aligned} & \lim_{\Delta t \to \infty} u_{i,j,k}(\Delta t, P = 2) = \\ & \lim_{\Delta t \to \infty} \left( e^{-2\bar{a}\Delta t} \left( u_{i,j,k}^n - \frac{\bar{d}}{\bar{a}} \left( \frac{a_0}{2} - \frac{a_1}{4\bar{a}} + \frac{a_2}{4\bar{a}^2} \right) \right) + \frac{\bar{d}}{\bar{a}} \left( \frac{a_0}{2} + \frac{a_1}{2} \left( \Delta t - \frac{1}{2\bar{a}} \right) + \frac{a_2}{2} \left( \frac{1}{2\bar{a}^2} - \frac{\Delta t}{\bar{a}} + (\Delta t)^2 \right) \right) \right) = \\ & \lim_{\Delta t \to \infty} \left( \frac{\bar{d}}{\bar{a}} \left( \frac{a_0}{2} + \frac{a_1}{2} \left( \Delta t - \frac{1}{2\bar{a}} \right) + \frac{a_2}{2} \left( \frac{1}{2\bar{a}^2} - \frac{\Delta t}{\bar{a}} + (\Delta t)^2 \right) \right) \right) = \\ & \lim_{\Delta t \to \infty} \left( \frac{\bar{d}}{\bar{a}} \left( \frac{\mathscr{U}^n}{2} + \frac{(-3\mathscr{U}^n + 4\mathscr{U}^{n+\frac{1}{2}} - \mathscr{U}^{n+1})}{2\Delta t} \left( \Delta t - \frac{1}{2\bar{a}} \right) + \frac{(2\mathscr{U}^n - 4\mathscr{U}^{n+\frac{1}{2}} + 2\mathscr{U}^{n+1})}{2(\Delta t)^2} \left( \frac{1}{2\bar{a}^2} - \frac{\Delta t}{\bar{a}} + (\Delta t)^2 \right) \right) \right) = \\ & \lim_{\Delta t \to \infty} \left( \frac{\bar{d}}{\bar{a}} \left( \frac{\mathscr{U}^n}{2} + \frac{(-3\mathscr{U}^n + 4\mathscr{U}^{n+\frac{1}{2}} - \mathscr{U}^{n+1})}{2} + \frac{(2\mathscr{U}^n - 4\mathscr{U}^{n+\frac{1}{2}} + 2\mathscr{U}^{n+1})}{2} \right) \right) = \frac{\bar{d}}{\bar{a}} \left( \frac{\mathscr{U}^{n+1}}{2} \right). \end{aligned} \tag{87}$$

For higher-order terms, i.e., for $P > 2$, we can verify that $\lim_{\Delta t \to \infty} u(\Delta t, P > 2) = \frac{\bar{d}}{\bar{a}} \cdot \frac{\mathscr{U}^{n+1}}{2}$.

Therefore, by induction, we conclude that $\lim_{\Delta t \to \infty} u(\Delta t, P = 0) = \frac{\bar{d}}{\bar{a}} \cdot \frac{\mathscr{U}^n}{2}$, and $\lim_{\Delta t \to \infty} u(\Delta t, P \neq 0) = \frac{\bar{d}}{\bar{a}} \cdot \frac{\mathscr{U}^{n+1}}{2}$. This substantiates the proposition. □

# 9 Temporal-ODE Discretization Scheme with Uncertain/Stochastic Diffusion Coefficient: $P^{th}$-Order Solution

We revisit the temporal-ODE scheme in the presence of an uncertain diffusion coefficient, denoted by $\tilde{\bar{d}}$. The $P^{th}$-order representative ODE retains nearly the same notation as in Eqs. (60–63), with the exception that we now define $\bar{\bar{a}} := \left(\frac{1}{(\Delta x)^2} + \frac{1}{(\Delta y)^2} + \frac{1}{(\Delta z)^2}\right)$, and the ODE becomes:

$$\frac{du}{d\tau} = \tilde{\bar{d}}\left(\sum_{p=0}^{P} a_p \tau^p - 2\bar{\bar{a}} u\right) + s. \tag{88}$$

Therefore, we can write the closed-form solution as:

$$u_{i,j,k}(\tau) = \left(u_{i,j,k}^n - \frac{s}{2\tilde{\bar{d}}\bar{\bar{a}}} - \frac{1}{\bar{\bar{a}}} \sum_{p=0}^{P} \frac{a_p}{2} \frac{(-1)^p p!}{(2\tilde{\bar{d}}\bar{\bar{a}})^p}\right) e^{-2\tilde{\bar{d}}\bar{\bar{a}}\tau} + \frac{1}{\bar{\bar{a}}} \sum_{p=0}^{P} \frac{a_p}{2} \sum_{q=0}^{p} \frac{\tau^{p-q}(-1)^q p!}{(2\tilde{\bar{d}}\bar{\bar{a}})^q (p-q)!} + \frac{s}{2\tilde{\bar{d}}\bar{\bar{a}}}. \tag{89}$$

where $\tau \in [0, \Delta t]$, and $u_{i,j,k}^n := u_{i,j,k}(0)$. For simplicity, we assume that $\tilde{\bar{d}}$ has a uniform probability density function (PDF) over the range $\tilde{\bar{d}} \in [\tilde{\bar{d}}_l, \tilde{\bar{d}}_u]$, where the subscripts $l$ and $u$ represent the lower and upper bounds, respectively. Therefore, we have $\int_{\tilde{\bar{d}}_l}^{\tilde{\bar{d}}_u} f(\tilde{\bar{d}}) d\tilde{\bar{d}} = 1$, where $f(\tilde{\bar{d}}) = \frac{1}{\tilde{\bar{d}}_u - \tilde{\bar{d}}_l}$, and the expected value $\mathbb{E}[u_i(\tau)]$ is given by: $\mathbb{E}[u_i(\tau)] = \int_{\tilde{\bar{d}}_l}^{\tilde{\bar{d}}_u} u_{i,j,k}(\tau) f(\tilde{\bar{d}}) d\tilde{\bar{d}}$.

Substituting $u_{i,j,k}(\tau)$ into the expected value formula and integrating, we obtain the update formula at $t_{n+1}$ (or $\tau = \Delta t$) as:

$$\mathbb{E}[u_i^{n+1}] = \\ \left(u_i^n - \frac{a_0}{2}\right) \frac{e^{-2\tilde{\bar{d}}_l \bar{\bar{a}} \Delta t} - e^{-2\tilde{\bar{d}}_u \bar{\bar{a}} \Delta t}}{2\bar{\bar{a}} \Delta t (\tilde{\bar{d}}_u - \tilde{\bar{d}}_l)} - \frac{s}{2\bar{\bar{a}}} \frac{Ei(-2\bar{\bar{a}}\Delta t \tilde{\bar{d}}_u) - Ei(-2\bar{\bar{a}}\Delta t \tilde{\bar{d}}_l)}{\tilde{\bar{d}}_u - \tilde{\bar{d}}_l} + \frac{a_0}{2} + \frac{s}{2\bar{\bar{a}}(\tilde{\bar{d}}_u - \tilde{\bar{d}}_l)} \ln\left(\frac{\tilde{\bar{d}}_u}{\tilde{\bar{d}}_l}\right) - \\ \frac{1}{\bar{\bar{a}}(\tilde{\bar{d}}_u - \tilde{\bar{d}}_l)} \sum_{p=1}^{P} \frac{a_p}{2} \sum_{q=0}^{p} \frac{(\Delta t)^{p-q}(-1)^q p!}{(2\bar{\bar{a}})^q (p-q)!} \frac{\tilde{\bar{d}}_u^{1-q} - \tilde{\bar{d}}_l^{1-q}}{(1-q)} - \frac{1}{\bar{\bar{a}}(\tilde{\bar{d}}_u - \tilde{\bar{d}}_l)} \sum_{p=1}^{P} \frac{a_p}{2} \frac{(-1)^p p!}{(2\bar{\bar{a}})^p} \int_{\tilde{\bar{d}}_l}^{\tilde{\bar{d}}_u} \frac{e^{-2\tilde{\bar{d}}\bar{\bar{a}}\Delta t}}{\tilde{\bar{d}}^p} d\tilde{\bar{d}}. \tag{90}$$

where $Ei(\cdot)$ is the exponential function. In addition, the last integral ($\int_{\tilde{\bar{d}}_l}^{\tilde{\bar{d}}_u} \frac{e^{-2\tilde{\bar{d}}\bar{\bar{a}}\Delta t}}{\tilde{\bar{d}}^p} d\tilde{\bar{d}}$) can either be computed numerically or analytically using the properties of incomplete gamma function.

# 10 Temporal-ODE Discretization Scheme for Advection Operator

The pure advection equation can be solved using various numerical discretization techniques, each with their own advantages and drawbacks. However, the characteristic method provides the clearest insight into the nature of the advection operator, transporting the initial condition $u^n$ along characteristic lines to new spatial positions. However, the existing numerical discretization techniques inherently deviate from reflecting this structural property of the advection operator in their update formulas.

The directional-ODE discretization approach, on the other hand, can offer an analytical update formula that can inherently preserve the structural property of the operator. To illustrate this, consider the standard 1D advection equation:

$$\frac{\partial u}{\partial t} + c\frac{\partial u}{\partial x} = 0. \tag{91}$$

where $c$ is a constant and the initial condition is given as $u(x,0) = u_0(x)$.

By differentiating with respect to $x$ and $t$, we obtain the following relations:

$$\frac{\partial^2 u}{\partial t^2} + c\frac{\partial^2 u}{\partial t \partial x} = 0, \quad \frac{\partial^2 u}{\partial x \partial t} + c\frac{\partial^2 u}{\partial x^2} = 0. \tag{92}$$

Recognizing that the mixed partial derivatives commute, i.e., $\frac{\partial^2 u}{\partial t \partial x} = \frac{\partial^2 u}{\partial x \partial t}$, one arrives at the following wave equation:

$$\frac{\partial^2 u}{\partial t^2} - c^2\frac{\partial^2 u}{\partial x^2} = 0, \quad u(x,0) = u_0(x), \quad \frac{\partial u}{\partial t}\bigg|_{t=0} = -c\frac{\partial u}{\partial x}\bigg|_{t=0}. \tag{93}$$

Exploiting the temporal-ODE scheme, we write:

$$\frac{d^2 u_i}{dt^2} = c^2 \frac{u_{i-1}(t) - 2u_i(t) + u_{i+1}(t)}{(\Delta x)^2}, \quad t \in [t_n, t_{n+1}]. \tag{94}$$

For clarity and simplicity, we consider the values of $u_{i-1}$ and $u_{i+1}$ at time $t_n$. Additionally, we define $\tau := t - t_n$, leading to the following ODE:

$$\frac{d^2 u_i}{d\tau^2} = c^2 \frac{u_{i-1}^n - 2u_i(\tau) + u_{i+1}^n}{(\Delta x)^2}, \quad \tau \in [0, \Delta t], \quad \Delta t := t_{n+1} - t_n. \tag{95}$$

The above equation can be written as:

$$\frac{d^2 u_i}{d\tau^2} = -Au_i(\tau) + B, \quad u_i(0) = u_i^n, \quad \frac{du_i}{d\tau}\bigg|_{\tau=0} = -c\frac{u_{i+1}^n - u_{i-1}^n}{2\Delta x}. \tag{96}$$

where $A = 2\frac{c^2}{(\Delta x)^2}$, and $B = \frac{u_{i-1}^n + u_{i+1}^n}{(\Delta x)^2}$.

The analytic solution to the above equation is:

$$u_i(\tau) = (u_i^n - \frac{B}{A})\cos(\theta) + \frac{1}{\sqrt{A}}\frac{du_i}{d\tau}\Big|_{\tau=0}\sin(\theta) + \frac{B}{A}, \quad \theta := \tau\sqrt{A}. \tag{97}$$

Therefore, the update formula at $t^{n+1}$ reads:

$$u_i^{n+1} = (u_i^n - \frac{B}{A})\cos(\theta) + \frac{1}{\sqrt{A}}\frac{du_i}{d\tau}\Big|_{\tau=0}\sin(\theta) + \frac{B}{A}, \quad \theta := \Delta t\sqrt{A}. \tag{98}$$

For numerical implementation, $\theta$ is typically chosen in the first quadrant.

**Remark:** For improved accuracy, a multi-stage predictor-corrector scheme can be easily derived for the above discrete scheme, where the corrector step approximates the neighboring points using a $P^{\text{th}}$-order polynomial.

It is evident that the above update formula reveals a propagating-wave structure of the operator.

# 11 The Directional-ODE Discretization without Splitting Techniques

The directional-ODE discretization approach can be applied to advection-diffusion equations without splitting techniques. However, as discussed throughout the paper, the representative ODEs can be derived in various ways, each with its own features. Here, to demonstrate this potential solely, we present some preliminary analysis using the temporal-ODE and spatial-ODE schemes.

## 11.1 The Case of a Linear ADE: Temporal-ODE Discretization Scheme

Let us begin by illustrating the method with the standard one-dimensional linear ADE:

$$\frac{\partial u}{\partial t} + c\frac{\partial u}{\partial x} = D\frac{\partial^2 u}{\partial x^2}. \tag{99}$$

where $c$ is a constant and the initial condition is given as $u(x,0) = u_0(x)$.

To incorporate both advection and diffusion operators into the discrete formulation while preserving the unconditional convergence property, we discretize the one-dimensional ADE analytically as:

**a)** $c > 0 \Rightarrow$ Backward Discretization for $\frac{\partial u}{\partial x}$:

$$\frac{du_i}{d\tau} = D\frac{u_{i-1}(\tau) + u_{i+1}(\tau)}{(\Delta x)^2} + c\frac{u_{i-1}(\tau)}{\Delta x} + u_i(\tau)\Big(-\frac{2D}{(\Delta x)^2} - \frac{c}{\Delta x}\Big), \quad \tau \in [0, \Delta t]. \tag{100}$$

**b)** $c < 0 \Rightarrow$ Forward Discretization for $\dfrac{\partial u}{\partial x}$ :

$$\frac{du_i}{d\tau} = D\frac{u_{i-1}(\tau) + u_{i+1}(\tau)}{(\Delta x)^2} - c\frac{u_{i+1}(\tau)}{\Delta x} + u_i(\tau)\left(-\frac{2D}{(\Delta x)^2} + \frac{c}{\Delta x}\right), \quad \tau \in [0, \Delta t]. \tag{101}$$

The above switching scheme ensures that the overall discrete form remains unconditionally stable. Specifically, considering the neighboring points at $t_n$, we can write the zeroth-order discrete schemes as the following linear ODEs:

$$\begin{aligned}\textbf{a) } c > 0: \quad &\frac{du_i}{d\tau} = Au_i(\tau) + B, \quad A = -\frac{2D}{(\Delta x)^2} - \frac{c}{\Delta x}, \quad B = D\frac{u_{i-1}^n + u_{i+1}^n}{(\Delta x)^2} + c\frac{u_{i-1}^n}{\Delta x}, \\ \textbf{b) } c < 0: \quad &\frac{du_i}{d\tau} = A'u_i(\tau) + B', \quad A' = -\frac{2D}{(\Delta x)^2} + \frac{c}{\Delta x}, \quad B' = D\frac{u_{i-1}^n + u_{i+1}^n}{(\Delta x)^2} - c\frac{u_{i+1}^n}{\Delta x}.\end{aligned} \tag{102}$$

The update formulas, representing the solutions at $\tau = \Delta t$, are given by:

$$\begin{aligned}\textbf{a) } c > 0: \quad &u_i^{n+1} = -\frac{B}{A} + \left(\frac{B}{A} + u_i^n\right)e^{A\Delta t}, \\ \textbf{b) } c < 0: \quad &u_i^{n+1} = -\frac{B'}{A'} + \left(\frac{B'}{A'} + u_i^n\right)e^{A'\Delta t}\end{aligned} \tag{103}$$

Since $A$, and $A'$ are both negative, the switching discrete scheme remains unconditionally stable.

**Remark:** A multi-stage predictor-corrector scheme can be readily derived (where the corrector step approximates the neighboring points using a $P^{th}$-order polynomial) for the above switching discrete scheme.

## 11.2 The Case of a Linear ADE: Spatial-ODE Discretization Scheme

Using the spatial-ODE scheme, the one-dimensional ADE can be written as:

$$\frac{u_i^{n+1} - u_i^n}{\Delta t} + c\frac{du_i}{dx}\Big|^{n+1} = D\frac{d^2 u_i}{dx^2}\Big|^{n+1}. \tag{104}$$

Therefore, the representative spatial ODE at $t_{n+1}$ can be written as:

$$\frac{d^2 u}{dx^2} - C\frac{du}{dx} - Au = B. \tag{105}$$

where, $A = \frac{1}{D\Delta t}$, $B = -\frac{u_i^n}{D\Delta t}$, and $C = \frac{c}{D}$.

In order to form a spatial-ODE scheme, we solve the above ODE on $\Omega := [x_i - \Delta x, x_i + \Delta x]$. In other words, the neighboring nodes are selected as the required boundary conditions:

$$u(x_i - \Delta x) = u_{i-1}^{n+1}, \quad u(x_i + \Delta x) = u_{i+1}^{n+1}. \tag{106}$$

The analytic solution for the above ODE reads:

$$u(x) = -\frac{B}{A} + c_1 e^{\lambda_1 x} + c_2 e^{\lambda_2 x}. \tag{107}$$

where,

$$\lambda_1 = \frac{C + \sqrt{C^2 + 4A}}{2}, \quad \lambda_2 = \frac{C - \sqrt{C^2 + 4A}}{2}. \tag{108}$$

Applying the boundary conditions, we can obtain the solution at $u(x_i) \equiv u_i^{n+1}$ as:

$$u_i^{n+1} = k_1 u_{i-1}^{n+1} + k_2 u_{i+1}^{n+1} + k_3. \tag{109}$$

where,

$$\begin{aligned} k_1 &= \frac{e^{\lambda_2 \Delta x} - e^{\lambda_1 \Delta x}}{e^{\Delta x(\lambda_2 - \lambda_1)} - e^{\Delta x(\lambda_1 - \lambda_2)}}, \quad k_2 = \frac{e^{-\lambda_1 \Delta x} - e^{-\lambda_2 \Delta x}}{e^{\Delta x(\lambda_2 - \lambda_1)} - e^{\Delta x(\lambda_1 - \lambda_2)}}, \\ k_3 &= -\frac{B}{A} + \frac{\frac{B}{A}\left(e^{\lambda_2 \Delta x} - e^{\lambda_1 \Delta x} + e^{-\lambda_1 \Delta x} - e^{-\lambda_2 \Delta x}\right)}{e^{\Delta x(\lambda_2 - \lambda_1)} - e^{\Delta x(\lambda_1 - \lambda_2)}} \end{aligned} \tag{110}$$

The above discrete formula can be solved using TDMA.

## 11.3 The Case of a Nonlinear ADE: Temporal-ODE Discretization Scheme

We begin by considering the following nonlinear one-dimensional ADE:

$$\frac{\partial u}{\partial t} + c(u)\frac{\partial u}{\partial x} = D\frac{\partial^2 u}{\partial x^2}, \tag{111}$$

where $c(u)$ is a nonlinear function of the field variable $u$.

In this context, the analysis developed for the linear case can be directly adapted by replacing $c$ with $c(u_i)$ evaluated at time step $t_n$. However, it is also possible to derive nonlinear representative ODEs that characterize the system dynamics more precisely. For brevity, we outline the form of the representative ODEs for three commonly used discretization schemes:

$$\begin{aligned} \text{Backward Discretization for } \frac{\partial u}{\partial x}: \quad & \frac{du}{d\tau} = Au + c(u)\left(Bu + C\right) + D, \\ \text{Forward Discretization for } \frac{\partial u}{\partial x}: \quad & \frac{du}{d\tau} = A'u + c(u)\left(B'u + C'\right) + D', \\ \text{Central Discretization for } \frac{\partial u}{\partial x}: \quad & \frac{du}{d\tau} = A''u + c(u)\left(B''\right) + C''. \end{aligned} \tag{112}$$

Here, the coefficients $A, B, C, D, A', B', C', D', A'', B'', C''$ may be constants or functions of the neighboring points, reflecting the nature of the discretization scheme employed.

**Remark 1:** While closed-form solutions are generally preferred due to their accuracy and the insights they provide into the behavior of the discrete scheme, obtaining closed-form solutions for the nonlinear case may not always be feasible. In such instances, alternative semi-analytical methods, such as the SADM, can serve as an effective approach to derive approximate analytical solutions.

**Remark 2:** As with previous temporal-ODE schemes, a multi-stage predictor-corrector approach can be designed to enhance the accuracy of the numerical solution. In this approach, the neighboring points are approximated using a $P^{\text{th}}$-order polynomial.

## 11.4 The Case of a Nonlinear ADE: Spatial-ODE Discretization Scheme

For the nonlinear case described by Eq. (111), the analysis for the linear case can again be readily extended by replacing $c$ with $c(u_i)$ at time step $t_n$. However, it is also possible to formulate a nonlinear representative ODE that more accurately reflects the system dynamics. The nonlinear representative ODE is expressed as:

$$\frac{d^2 u}{dx^2} + Ac(u)\frac{du}{dx} + Bu = C. \tag{113}$$

where $A, B, C$ are either constants or functions of the previous time step $t_n$.

Notably, the nonlinear nature of the term $c(u)$ in this equation often precludes the derivation of a closed-form solution. In such cases, semi-analytical methods, such as SADM, offer a practical alternative to obtain simple yet flexible analytic approximations, which can vary in order and provide insights into the solution's behavior under different conditions.

## 12 ADI Method for $\psi - \omega$ Navier Stokes Equations

The Navier-Stokes equations in the streamfunction-vorticity ($\psi$-$\omega$) formulation for 2D incompressible flow are given as follows:

$$\begin{aligned}
\frac{\partial \omega}{\partial t} + \frac{\partial \psi}{\partial y}\frac{\partial \omega}{\partial x} - \frac{\partial \psi}{\partial x}\frac{\partial \omega}{\partial y} &= \frac{1}{Re}\left(\frac{\partial^2 \omega}{\partial x^2} + \frac{\partial^2 \omega}{\partial y^2}\right), \\
\frac{\partial^2 \psi}{\partial x^2} + \frac{\partial^2 \psi}{\partial y^2} &= -\omega.
\end{aligned} \tag{114}$$

where $\omega$ denotes the vorticity, defined as $\omega = \frac{\partial v}{\partial x} - \frac{\partial u}{\partial y}$, and $Re$ represents the Reynolds number, expressed as $Re = \frac{UW}{\nu}$, where $U$ is the characteristic velocity, $W$ is the characteristic length (which is the channel width), and $\nu$ is the kinematic viscosity. The streamfunction $\psi$ is related to the velocity components through the relationships $u = \frac{\partial \psi}{\partial y}$, and $v = -\frac{\partial \psi}{\partial x}$ where $u$ and $v$ are the normalized velocity components in the $x$ and $y$ directions, respectively.

Here, we provide a concise summary of the ADI algorithm employed to solve the $\psi$-$\omega$ Navier-Stokes equations.

**x-direction discretization for $\omega$:**

$$\frac{\omega_{i,j}^* - \omega_{i,j}^n}{\frac{\Delta t}{2}} + u^n \frac{\omega_{i+1,j}^* - \omega_{i,j}^*}{2\Delta x} + v^n \frac{\omega_{i,j+1}^n - \omega_{i,j-1}^n}{2\Delta y} = \\ \frac{1}{Re}\Big(\frac{\omega_{i+1,j}^* - 2\omega_{i,j}^* + \omega_{i-1,j}^*}{(\Delta x)^2} + \frac{\omega_{i,j+1}^n - 2\omega_{i,j}^n + \omega_{i,j-1}^n}{(\Delta y)^2}\Big). \tag{115}$$

**y-direction discretization for $\omega$:**

$$\frac{\omega_{i,j}^{n+1} - \omega_{i,j}^*}{\frac{\Delta t}{2}} + u^n \frac{\omega_{i+1,j}^* - \omega_{i,j}^*}{2\Delta x} + v^n \frac{\omega_{i,j+1}^{n+1} - \omega_{i,j-1}^{n+1}}{2\Delta y} = \\ \frac{1}{Re}\Big(\frac{\omega_{i+1,j}^* - 2\omega_{i,j}^* + \omega_{i-1,j}^*}{(\Delta x)^2} + \frac{\omega_{i,j+1}^{n+1} - 2\omega_{i,j}^{n+1} + \omega_{i,j-1}^{n+1}}{(\Delta y)^2}\Big). \tag{116}$$

**Guess-Seidel discretization for $\psi$**

$$\frac{\psi_{i+1,j}^n - 2\psi_{i,j}^{n+1} + \psi_{i-1,j}^n}{(\Delta x)^2} + \frac{\psi_{i,j+1}^n - 2\psi_{i,j}^{n+1} + \psi_{i,j-1}^n}{(\Delta y)^2} = -\omega_{i,j}^{n+1}. \tag{117}$$

Therefore, the numerical algorithm involves:

**TDMA in x-direction, 1 iteration:**

$$\omega_{i-1,j}^*\Big(\lambda_1 + \frac{u_{i,j}^n \Delta t}{4\Delta x}\Big) + \omega_{i,j}^*\big(-1 - 2\lambda_1\big) + \omega_{i+1,j}^*\Big(\lambda_1 - \frac{u_{i,j}^n \Delta t}{4\Delta x}\Big) = \\ -\omega_{i,j}^n + v_{i,j}^n \frac{\Delta t}{2}\frac{\omega_{i,j+1}^n - \omega_{i,j-1}^n}{2\Delta y} - \frac{\Delta t}{2Re}\frac{\omega_{i,j+1}^n - 2\omega_{i,j}^n + \omega_{i,j-1}^n}{(\Delta y)^2}. \tag{118}$$

**TDMA in y-direction, updating $\omega$, 1 iteration:**

$$\omega_{i,j-1}^{n+1}\Big(\lambda_2 + \frac{v_{i,j}^n \Delta t}{4\Delta x}\Big) + \omega_{i,j}^{n+1}\big(-1 - 2\lambda_2\big) + \omega_{i,j-1}^{n+1}\Big(\lambda_2 - \frac{v_{i,j}^n \Delta t}{4\Delta x}\Big) = \\ -\omega_{i,j}^* + u_{i,j}^n \frac{\Delta t}{2}\frac{\omega_{i+1,j}^* - \omega_{i-1,j}^*}{2\Delta x} - \frac{\Delta t}{2Re}\frac{\omega_{i+1,j}^* - 2\omega_{i,j}^* + \omega_{i-1,j}^*}{(\Delta x)^2}. \tag{119}$$

**Updating $\psi$, $K$ iterations, $K \geq 1$:**

$$\psi_{i,j}^{n+1} = \frac{\Big(\frac{\psi_{i+1,j}^n + \psi_{i-1,j}^n}{(\Delta x)^2} + \frac{\psi_{i,j+1}^n + \psi_{i,j-1}^n}{(\Delta y)^2} + \omega_{i,j}^{n+1}\Big)}{\frac{2}{(\Delta x)^2} + \frac{2}{(\Delta y)^2}}. \tag{120}$$

The implementation of the ADI scheme described above is available in link.